\documentclass[leqno,11pt]{article}\textwidth 160mm\textheight 235mm
\usepackage{amsfonts}
\usepackage{amssymb}
\usepackage{graphicx}
\usepackage{amsmath}
\usepackage{amsthm}
\usepackage{enumitem}
\usepackage{mathrsfs}
\usepackage{tikz}
\usepackage{cancel}
\usepackage{todonotes}
\usetikzlibrary{matrix}
\usetikzlibrary{arrows,shapes}
\usepackage{cite}
\usepackage{hyperref}
\hypersetup{colorlinks=true,urlcolor=blue}
\expandafter\let\expandafter\oldproof\csname\string\proof\endcsname
\let\oldendproof\endproof
\renewenvironment{proof}[1][\proofname]{
\oldproof[\ttfamily\scshape \bf #1.]
}{\oldendproof}
\def\ve{\varepsilon}

\def\tilde{\widetilde}
\def\emp{\emptyset}

\def\dom{{\rm dom}\,}
\def\span{{\rm span}\,}
\def\epi{{\rm epi\,}}
\def\bd{{\rm bd\,}}

\def\min{\mbox{\rm minimize}}

\def\Lm{{\Lambda}}
\def\d{{\rm d}}
\def\sub{\partial}
\def \Q{{\cal Q}}
\def\B{\mathbb B}
\def\L{{\mathscr{L}}}

\def\ox{\overline{x}}
\def\oy{\overline{y}}
\def\oz{\overline{z}}

\def\disp{\displaystyle}
\def\tto{\rightrightarrows}
\def\Hat{\widehat}

\def\Bar{\overline}
\def\ra{\rangle}
\def\la{\langle}
\def\ve{\varepsilon}
\def\epsilon{\varepsilon}
\def\ox{\bar{x}}
\def\oy{\bar{y}}
\def\oz{\bar{z}}

\def\ov{\bar{v}}

\def\ou{\bar{u}}

\def\co{\mbox{\rm co}\,}

\def\inte{\mbox{\rm int}\,}
\def\gph{\mbox{\rm gph}\,}
\def\epi{\mbox{\rm epi}\,}

\def\dom{\mbox{\rm dom}\,}
\def\ker{\mbox{\rm ker}\,}

\def\dn{\downarrow}
\def\O{\Omega}
\def\ph{\varphi}
\def\emp{\emptyset}
\def\st{\stackrel}
\def\oR{\Bar{\R}}
\def\lm{\lambda}
\def\gg{\gamma}
\def\dd{\delta}
\def\al{\alpha}
\def\Th{\Theta}
\def \N{{\rm I\!N}}
\def \R{{\rm I\!R}}

\def\vep{\varepsilon}

\def\Limsup{\mathop{{\rm Lim}\,{\rm sup}}}
\def\sm{\hbox{${1\over 2}$}}

\def\sce{\setcounter{equation}{0}}
\begin{document}
\newtheorem{Theorem}{Theorem}[section]
\newtheorem{Proposition}[Theorem]{Proposition}
\newtheorem{Remark}[Theorem]{Remark}
\newtheorem{Lemma}[Theorem]{Lemma}
\newtheorem{Corollary}[Theorem]{Corollary}
\newtheorem{Definition}[Theorem]{Definition}
\newtheorem{Example}[Theorem]{Example}
\newtheorem{Algorithm}[Theorem]{Algorithm}
\renewcommand{\theequation}{{\thesection}.\arabic{equation}}
\renewcommand{\thefootnote}{\fnsymbol{footnote}}
\begin{center}
{\bf\Large Generalized Newton Algorithms for Tilt-Stable Minimizers\\ in Nonsmooth Optimization}\\[1ex]
BORIS S. MORDUKHOVICH\footnote{Department of Mathematics, Wayne State University, Detroit, MI 48202, USA (boris@math.wayne.edu). Research of this author was partly supported by the US National Science Foundation under grants DMS-1512846 and DMS-1808978, by the US Air Force Office of Scientific Research grant \#15RT04, and by the Australian Research Council under Discovery Project DP-190100555.} and M. EBRAHIM SARABI\footnote{Department of Mathematics, Miami University, Oxford, OH 45065, USA (sarabim@miamioh.edu).}
\end{center}
\small{\bf Abstract.} This paper aims at developing two versions of the generalized Newton method to compute local minimizers for nonsmooth problems of unconstrained and constrained optimization that satisfy an important stability property known as tilt stability. We start with unconstrained minimization of continuously differentiable cost functions having Lipschitzian gradients and suggest two second-order algorithms of the Newton type: one involving coderivatives of Lipschitzian gradient mappings, and the other based on graphical derivatives of the latter. Then we proceed with the propagation of these algorithms to minimization of extended-real-valued prox-regular functions, while covering in this way problems of constrained optimization, by using Moreau envelopes. Employing advanced techniques of second-order variational analysis and characterizations of tilt stability allows us to establish the solvability of subproblems in both algorithms and to prove the $Q$-superlinear convergence of their iterations.\\[1ex]
{\bf Key words.} Nonsmooth optimization, generalized Newton method, tilt-stable local minimizers, prox-regular functions, superlinear convergence.\\[1ex]
{\bf  Mathematics Subject Classification (2000)} 90C31, 49J52, 49J53\\[1ex]
{\bf Abbreviated Title}. Newton Method for Tilt-Stable Minimizers\vspace*{-0.15in}

\normalsize
\section{Introduction}\label{intro}\sce\vspace*{-0.05in}

The classical Newton methods for solving equations and optimization problem, as well as their various modifications and extensions, have been well recognized among the most efficient numerical algorithms to find local solutions; see, e.g., the books \cite{dr,fp,is14,kk} with the vast commentaries and references therein. The standard framework of the Newton methods, which goes back to the Newton method of tangents, is to solve smooth equations $g(x)=0$ with $g\colon\R^n\to\R^n$. It is then applied to finding local solutions for problems of unconstrained optimization of the type
\begin{equation}\label{main}
\mbox{ minimize }\;\ph(x)\;\mbox{ subject to }\;x\in\R^n
\end{equation}
with ${\cal C}^2$-smooth objective/cost functions $\ph\colon\R^n\to\R$ by solving the stationary equations $g(x):=\nabla\ph(x)=0$ based on the classical Fermat necessary optimality condition. The corresponding                                                                                                                                                                                                                                                                                                Newton algorithm designed in this way is expressed via the Hessian matrix of $f$ at a solution point and exhibits local superlinear convergence when the Hessian matrix $\nabla^2\ph(\ox)$ at a solution point is positive-definite.

Among various extensions of the Newton algorithm to solve nonsmooth equations $g(x)=0$ with Lipschitzian mappings, the most successful one is the so-called {\em semismooth Newton method} initiated independently by Kummer \cite{ku88} and by Qi and Sun \cite{qs93}. Applying to problems of unconstrained optimization \eqref{main}, the semismooth Newton method addresses objective functions of class ${\cal C}^{1,1}$ (labeled also as ${\cal C}^{1+})$ around local minimizers $\ox$, i.e., the class of ${\cal C}^1$-smooth functions with locally Lipschitzian gradients. For this important class of nonsmooth (of the second-order) problems, superlinear convergence of the semismooth Newton method was achieved under some additional requirements; see \cite{fp,is14,kk} and the discussions below.

An interesting idea to extend the semismooth Newton method to problems of convex optimization without the ${\cal C}^{1,1}$ requirement on cost functions was suggested by Fukushima and Qi \cite{fq} in the framework of \eqref{main} with a nondifferentiable convex objective $\ph\colon\R^n\to\R$ finite on the whole space. They proposed to consider a regularization of problem \eqref{main} by replacing $\ph$ with its Moreau envelope. It is well known in convex analysis that the latter function is real-valued, convex, and everywhere differentiable on $\R^n$. Employing in this way the machinery of the semismooth Newton method, superlinear convergence of the corresponding algorithm was achieved in \cite{fq} under appropriate regularity assumptions.\vspace*{0.03in}

In this paper we offer a novel viewpoint on developing Newton methods in both unconstrained and constrained optimization and design new Newton-type algorithms for nonsmooth optimization problems. To the best of our knowledge, for the first time in the literature we aim at designing algorithms that seek not roots of equations, but specifically concern optimization problems and address fast convergence to {\em stable} local minimizers. Our stability choice is the concept of {\em tilt stability} introduced by Poliquin and Rockafellar \cite{pr98} in the general extended-real-valued framework of unconstrained optimization, which implicitly incorporates constraints via effective domains of cost functions. It is shown in \cite{pr98} that for ${\cal C}^2$-smooth functions $\ph$ in \eqref{main}, any root/stationary point of $\nabla\ph(\ox)=0$ is a tilt-stable local minimizer of \eqref{main} if and only if the Hessian matrix $\nabla^2\ph(\ox)$ is positive-definite. This confirms that the positive-definiteness of $\nabla^2\ph(\ox)$, which ensures superlinear convergence of the classical Newton method, automatically brings us to tilt-stable local minimizers. In the more general cases investigated in this paper, we explicitly impose the tilt stability requirement in our Newton schemes developed below.

Independently of algorithmic applications, an important advantage of tilt-stable local minimizers is their currently achieved comprehensive {\em second-order characterizations} in both unconstrained and constrained formats of optimization. The first characterization of tilt stability was obtained in the seminal paper by Poliquin and Rockafellar \cite{pr98} for the class of prox-regular and subdifferentially continuous extended-real-valued functions $\ph\colon\R^n\to\oR:=(-\infty,\infty]$ introduced by them a bit earlier \cite{pr96}; see the exact definitions of these and other needed notions in Section~\ref{sect02}. It is the major class of functions used in second-order variational analysis. The tilt stability characterization obtained in \cite{pr98} was expressed precisely at the reference local minimum point in terms of the positive-definiteness of the {\em second-order subdifferential/generalized Hessian} of $\ph$ in the sense of Mordukhovich \cite{m92} defined via the coderivative of his first-order limiting subdifferential. Subsequent second-order characterizations of tilt stability and related notions have been established more recently for various classes of unconstrained and constrained optimization problems; see \cite{bgm,chn,dl,dmn,gm15,mn14,mnr,mos,ms14} among other publications. Some of these characterizations will be used in what follows for the design and justification of the suggested Newton-type algorithms that superlinearly converge to tilt-stable local minimizers.\vspace*{0.03in}

To proceed in this direction, we begin with optimization problems \eqref{main}, where the objective functions $\ph$ are of class ${\cal C}^{1,1}$ around the reference points. For such problems of unconstrained optimization, we design two independent second-order algorithms of the Newton type. The first algorithm is based on using the {\em coderivative} of the gradient mapping for $\ph$ (i.e., the aforementioned second-order subdifferential), while the second one employs the {\em graphical derivative} of $\nabla \ph$. We justify the {\em solvability} of subproblems in both algorithms for tilt-stable minimizers (with an additional twice epi-differentiability assumption on $\ph$ needed for the efficient realization of the second algorithm) and achieve their local {\em superlinear convergence} under the {\em semismoothness$^*$} of $\nabla\ph$, a property that has been recently introduced by Gfrerer and Outrata \cite{go19}. Note that the graphical derivative has been already used (from different prospectives) in generalized Newton methods to solve equations and inclusions in the general scheme of \cite{kk} and specifically in \cite{ds,hkmp}, but the coderivative-based algorithm seems to be completely novel in numerical optimization.

Next we turn, for the first time in the literature, to developing Newton-type algorithms to find tilt-stable minimizers of extended-real-valued {\em prox-regular} and {\em subdifferentially continuous} functions $\ph\colon\R^n\to\oR$ in \eqref{main} while encompassing in this way problems of constrained optimization, which are also considered explicitly in what follows. The main idea here is to reduce such problems to those with ${\cal C}^{1,1}$ objectives by using {\em Moreau envelopes}. To justify this procedure, we show that the tilt stability of local minimizers of $\ph$ and the semismoothness$^*$
of the subgradient mapping $\partial\ph$ in \eqref{main} generates the corresponding properties of the regularized one via the Moreau envelope of $\ph$. This allows us to establish the solvability of subproblems and superlinear convergence of both coderivative-based and graphical derivative-based generalized Newton algorithms for tilt-stable local minimizers of the major class of extended-real-valued functions under consideration.

Finally, we present applications of the above results obtained in the extended-real-valued format of unconstrained optimization to the class of explicitly constrained optimization problems written in the form of {\em conic programming} (although the underlying set $\Theta$ may not be a cone):
\begin{equation}\label{coop}
\mbox{minimize }\;\psi(x)\;\mbox{ subject to }\;f(x)\in\Th,
\end{equation}
where $\psi\colon\R^n\to\R$ and $f\colon\R^n\to\R^m$ are ${\cal C}^2$-smooth, and where $\Th\subset\R^m$ is closed and convex. Applying to such problems, the developed graphical derivative-based algorithm generates a new Newton-type algorithm involving {\em second subderivatives} of cost functions. Its justification and convergence analysis for tilt-stable minimizers employ the recent developments of \cite{mms19} on {\em parabolic regularity} in second-order variational analysis.\vspace*{0.03in}

It is important to emphasize that the generalized differential constructions used in the proposed algorithms  employ not just coderivatives and graphical derivatives for general set-valued mappings, but those applied to the classical gradient mappings. As it has been recently realized, such second-order constructions enjoy excellent calculus rules as {\em equalities} with efficient computation; see, e.g., \cite{mms19a,mms19,m18,mr} and the references therein. This makes the proposed algorithms more attractive for numerical implementations.\vspace*{0.03in}

The rest of the paper is organized as follows. Section~\ref{sect02} recalls major concepts of variational analysis used in this paper and presents some related preliminary results. In Section~\ref{sect03} we derive basic estimates needed below to verify the performance of the suggested Newton-type algorithms. Section~\ref{sect04} is devoted to the coderivative-based generalized Newton algorithm for ${\cal C}^{1,1}$ functions with verifying the solvability of its subproblems and superlinear convergence of its iterates to a tilt-stable minimizer. In Section~\ref{sect04a} we do the same for a Newton algorithm dealing with ${\cal C}^{1,1}$ functions that is based on graphical derivatives. Section~\ref{sect05} develops Newton algorithms of both type for the class of extended-real-valued prox-regula functions. The final Section~\ref{sect05a} provides the applications of the developed results to optimization problems with explicit constraints.\vspace*{0.03in}

Throughout this paper we employ standard notation of variational analysis and optimization; see, e.g., \cite{fp,m18,rw}. Recall that $\B$ stands for the closed unit ball in the space in question, $\B_r(x):=x+r\B$ is the closed ball centered at $x$ with radius $r>0$, and $\N:=\{1,2,\ldots\}$. For the reader's convenience and notational unification we use as a rule small Greek letters to denote scalar and extended-real-valued functions, small Latin letters for vectors and single-valued mappings, and capital letters for sets, set-valued mappings, and matrices. We also distinguish in notation between single-valued mappings $f\colon\R^n\to\R^m$ and set-valued ones $F\colon\R^n\tto\R^m$. The (Painlev\'e-Kuratowski) {\em outer limit} of a $F\colon\R^n\tto\R^m$ as $x\to\ox$ is defined by
\begin{equation}\label{pk}
\Limsup_{x\to\ox}F(x):=\big\{y\in\R^n\big|\;\exists\,x_k\to\ox,\;y_k\to y\;\mbox{ with }\;y_k\in F(x_k),\;k\in\N\big\}.
\end{equation}
Given a set $\O\subset\R^n$, its {\em indicator function} is defined by $\dd_\O(x):=0$ for $x\in\O$ and $\dd_\O(x):=\infty$ otherwise, while the {\em distance} from $x$ to $\O$ is denoted by dist$(x;\O)$. For a function $\ph\colon\R^n\to\R$, denote by $\nabla\ph(\ox)$ and $\nabla^2\ph(\ox)$ its gradient and Hessian at $\ox$, respectively. If $f=(f_1,\ldots,f_m)\colon\R^n\to\R^m$ is twice differentiable at $\ox\in\R^n$, its second derivative at $\ox$, labeled by $\nabla^2 f(\ox)$, is a bilinear mapping from $\R^n\times\R^n$ into $\R^m$ given by the representation
\begin{equation*}
\nabla^2f(\ox)(w,v):=\big(\la\nabla^2f_1(\ox)w,v\ra,\ldots,\la\nabla^2f_m(\ox)w,v\ra\big)\;\mbox{ for all }\;v,w\in\R^n.
\end{equation*}
\vspace*{-0.4in}

\section{Major Definitions and Preliminaries}\sce \label{sect02}\vspace*{-0.05in}

We begin by recalling some of well-known tools of variational analysis and generalized differentiation that will be utilized throughout this paper; see, e.g., \cite{m18,rw} for this and additional material. Given a nonempty set $\O\subset\R^n$ with $\ox\in\O$, the (Bouligand-Severi) {\em tangent/contingent cone} $T_ \O(\ox)$ to $\O$ at $\ox\in\O$ is defined by
\begin{equation}\label{tan}
T_\O(\ox):=\big\{w\in\R^n\big|\;\exists\,t_k{\dn}0,\;\;w_k\to w\;\;\mbox{ as }\;k\to\infty\;\;\mbox{with}\;\;\ox+t_kw_k\in\O\big\}.
\end{equation}
We say that a tangent vector $w\in T_\O(\ox)$ is {\em derivable} if there exists $\xi\colon[0,\ve]\to\O$ with $\ve>0$, $\xi(0)=\ox$, and $\xi'_+(0)=w$, where $\xi'_+$ stands for the right derivative of $\xi$
at $0$ given by
\begin{equation*}
\xi'_+(0):=\lim_{t\dn 0}\frac{\xi(t)-\xi(0)}{t}.
\end{equation*}

The (Fr\'echet) {\em regular normal cone} to $\O$ at $\ox\in\O$ is given by
\begin{equation}\label{rnc}
\Hat N_\O(\ox):=\disp\Big\{v\in\R^n\Big|\;\limsup_{x\st{\O}{\to}\ox}\frac{\la v,x-\ox\ra}{\|x-\ox\|}\le 0\Big\},
\end{equation}
where $x\st{\O}{\to}\ox$ indicates that $x\to\ox$ with $x\in\O$. The (Mordukhovich) {\em limiting normal cone} to the set $\O$ at $\ox\in\O$ is defined as the outer limit \eqref{pk} of \eqref{rnc} as $x\st{\O}{\to}\ox$ by
\begin{equation}\label{lnc}
N_\O(\ox):=\Limsup_{x\st{\O}{\to}\ox}\Hat N_\O(x).
\end{equation}

Given further an extended-real-valued function $\ph\colon\R^n\to\oR$ with
\begin{equation*}
\dom\ph:=\big\{x\in\R^n\big|\;\ph(x)<\infty\big\}\;\mbox{ and }\;\epi\ph:=\big\{(x,\al)\in\R^{n+1}\big|\;\al\ge\ph(x)\},
\end{equation*}
its (first-order limiting) {\em subdifferential} at $\ox\in\dom\ph$ is defined via \eqref{lnc} by
\begin{eqnarray}\label{1sub}
\partial\ph(\ox):=\big\{v\in\R^n\big|\;(v,-1)\in N_{{\scriptsize\epi\ph}}\big(\ox,\ph(\ox)\big)\big\}.
\end{eqnarray}
Recall also that the function $\ph$ is {\em proper} if $\dom\ph\ne\emp$.\vspace*{0.03in}

Considering further a set-valued mapping/multifunction $F\colon\R^n\tto\R^m$ with the graph $\gph F:=\{(x,y)\in\R^n\times\R^m|\;y\in F(x)\}$, its {\em graphical derivative} at $(\ox,\oy)\in\gph F$ is defined by using the tangent cone \eqref{tan} to the graph as
\begin{equation}\label{gder}
DF(\ox,\oy)(u):=\big\{v\in\R^m\big|\;(u,v)\in T_{\scriptsize{\gph F}}(\ox,\oy)\big\},\quad u\in\R^n.
\end{equation}
The {\em regular coderivative} and the {\em limiting coderivative} of $F$ at $(\ox,\oy)\in\gph F$ are defined via the normal cones \eqref{rnc} and \eqref{lnc}, respectively, by
\begin{equation}\label{rcod}
\Hat D^*F(\ox,\oy)(v):=\big\{u\in\R^n\big|\;(u,-v)\in\Hat N_{\scriptsize{\gph F}}(\ox,\oy)\big\},\quad v\in\R^m,
\end{equation}
\begin{equation}\label{lcod}
D^*F(\ox,\oy)(v):=\big\{u\in\R^n\big|\;(u,-v)\in N_{\scriptsize{\gph F}}(\ox,\oy)\big\},\quad v\in\R^m.
\end{equation}
When $F=f$ is single-valued, we drop $\oy=f(\ox)$ from the notation of the graphical derivative \eqref{gder} and both coderivatives in \eqref{rcod} and \eqref{lcod}.

Note that the limiting normal cone \eqref{lnc}, together with the subdifferential \eqref{1sub} and the coderivative \eqref{lcod} generated by it, enjoys comprehensive calculus rules based on variational and extremal principles of variational analysis; see the books \cite{m06,m18,rw} and the references therein. It is not the case for the regular normal cone \eqref{rnc} and the tangent cone \eqref{tan} as well as the generated regular coderivative \eqref{rcod} and graphical derivative \eqref{gder} constructions. However, quite recently it has been realized that the latter derivative construction, being applied to the first-order subdifferential mappings \eqref{1sub}, possesses nice calculation formulas in many important situations; see \cite{chn,gm15,go17,go19,mms19a,mms19,m18,mn14,mnr,mos,ms14} among other publications. In this way we enter the realm of {\em second-order} variational analysis and generalized differentiation.\vspace*{0.03in}

A major class of extended-real-valued functions, which overwhelmingly appears in second-order variational analysis and optimization, consists of prox-regular and subdifferentially continuous ones introduced in \cite{pr96}. Besides convex and ${\cal C}^2$-smooth functions, this family contains those which are strongly amenable, lower-${\cal C}^2$, of the maximum type, of class ${\cal C}^{1,1}$, etc.; see \cite{rw} for more details. Recall that $\ph\colon\R^n\to\oR$ is {\em prox-regular} at $\ox$ for $\ov$ if $\ph$ is finite at $\ox$ and locally lower semicontinuous (l.s.c.) around $\ox$ with $\ov\in\sub\ph(\ox)$, and there exist constants $\ve>0$ and $\rho\ge 0$ such that for all $x\in\B_{\vep}(\ox)$ with $\ph(x)\le\ph(\ox)+\ve$ we have
\begin{eqnarray}\label{prox}
\ph(x)\ge\ph(u)+\la\ov,x-u\ra-\frac{\rho}{2}\|x-u\|^2\;\mbox{ whenever }\;(u,v)\in(\gph\sub\ph)\cap\B_{\vep}(\ox,\ov).
\end{eqnarray}
The function $\ph$ is {\em subdifferentially continuous} at $\ox$ for $\ov$ if the convergence $(x_k,v_k)\to(\ox,\ov)$ with $v_k\in\sub\ph(x_k)$ yields $\ph(x_k)\to\ph(\ox)$ as $k\to\infty$. For brevity we say that $\ph$ is {\em continuously prox-regular} at $\ox$ for $\ov$ if it has both prox-regularity and subdifferential continuity properties. In this case, the condition $\ph(x)\le\ph(\ox)+\ve$ in the definition of prox-regularity can be omitted.

Next we define, following \cite{pr98}, the underlying notion of tilt stability for extended-real-valued functions. Given $\ph\colon\R^n\to\oR$, a point $\ox\in\dom\ph$ is said to be a {\em tilt-stable local minimizer} of the function $\ph$ if for some $\gg>0$ the argminimum mapping
\begin{equation}\label{tilt}
M_{\gg}\colon v\mapsto{\mbox{argmin}}\big\{\ph(x)-\la v,x\ra\big|\;x\in\B_{\gg}(\ox)\big\}
\end{equation}
is single-valued and Lipschitz continuous on a neighborhood of $\ov=0$ with $M_{\gg}(\ov)=\{\ox\}$. As in \cite{mn14}, we say that $\ox$ is a tilt-stable local minimizer for $f$
with {\em modulus} $\kappa>0$ if the mapping $M_\gg$ from \eqref{tilt} is Lipschitz continuous with constant $\kappa$ on a neighborhood of $\ov=0$ with $M_{\gg}(\ov)=\{\ox\}$.

Recall further that a set-valued mapping $F\colon\R^n\tto\R^m$ admits a {\em single-valued graphical localization} around $(\ox,\oy)\in\gph F$ if there exist some neighborhoods $U$ of $\ox$ and $V$ of $\oy$ together with a single-valued mapping $f\colon U\to V$ such that $(\gph F)\cap(U\times V)=\gph f$. Now we present a useful characterization of tilt stability taken from \cite[Theorem~3.2]{mn14}.\vspace*{-0.07in}

\begin{Proposition}[\bf tilt stability via the second-order growth condition]\label{usogc2} Let $\ph\colon\R^n\to\oR$ be continuously prox-regular at $\ox$ for $\ov=0$. The following are equivalent:

{\bf(i)} The point $\ox$ is a tilt-stable minimizer of $\ph$ with modulus $\kappa>0$.

{\bf(ii)} There are neighborhoods $U$ of $\ox$ and $V$ of $\ov$ such that the inverse mapping $(\sub\ph)^{-1}$ admits a single-valued localization $\vartheta\colon V\to U$ around $(\ov,\ox)$, and that for any pair $(v,u)\in\gph\vartheta=(\gph(\sub\ph^{-1}))\cap(V\times U)$ we have the {\sc uniform second-order growth condition}
\begin{equation}\label{usogc}
\ph(x)\geq\ph(u)+\la v,x-u\ra+\frac{1}{2\kappa}\|x-u\|^2\;\mbox{ whenever }\;x\in U.
\end{equation}
\end{Proposition}

Another useful characterization of tilt stability, taken from \cite[Theorem~2.1]{chn}, employs the graphical derivative \eqref{gder} of the subgradient mapping $\partial\ph$, which is a second-order generalized differential construction of variational analysis.\vspace*{-0.07in}

\begin{Proposition}[\bf tilt stability via the subgradient graphical derivative]\label{tiltch} Let $\ph\colon\R^n\to\oR$ be continuously prox-regular at $\ox$ for $\ov=0$. The following are equivalent:

{\bf(i)} The point $\ox$ is a tilt-stable local minimizer of $\ph$ with modulus $\kappa>0$.

{\bf(ii)} There exists a constant $\eta>0$ such that for all $w\in \R^n$ we have
\begin{equation*}
\la u,w\ra\ge\frac{1}{\kappa}\|w\|^2\;\mbox{ whenever }\;u\in\big(D\partial\ph\big)(x,v)(w)\;\mbox{ with }\;(x,v)\in(\gph\sub\ph)\cap\B_\eta(\ox,0).
\end{equation*}
\end{Proposition}

Next we consider functions $\ph\colon\R^n\to\R$ of class ${\cal C}^{1,1}$ around $\ox$, which are continuously prox-regular therein as defined above. The following proposition, taken from \cite[Theorem~4.7]{m06}, is used in the study of such functions via coderivatives.\vspace*{-0.07in}

\begin{Proposition}[\bf coderivatives of Lipschitzian mappings]\label{lipcon} Let $f\colon\R^n\to\R^m$ be Lipschitz continuous around $\ox\in\R^n$. Then there are positive numbers $\rho$ and $\eta$ such that
\begin{equation}\label{lip}
\|w\|\le\rho\|u\|\;\mbox{ for all }\;(u,w)\in\gph D^*f(x)\;\mbox{ and }\;x\in\B_\eta(\ox).
\end{equation}
\end{Proposition}

Continuing with single-valued mappings $f\colon\R^n\to\R^m$ that are locally Lipschitzian around $\ox$, define the collection of limiting Jacobian matrices
\begin{equation}\label{radem}
\Bar\nabla f(\ox):=\Big\{\lim_{k\to\infty}\nabla f(x_k)\Big|\;x_k\to\ox,\;x_k\in\O_f\Big\},
\end{equation}
where $\O_f$ stands for the set on which $f$ is differentiable. The classical Rademacher theorem tells us that \eqref{radem} is a nonempty compact in $\R^{m\times n}$. The (Clarke) {\em generalized Jacobian} of $f$ at $\ox$ is defined as the convex hull of the limiting Jacobian set \eqref{radem} and is denoted by $\co\Bar\nabla f(\ox)$. The following relationship is well known:
\begin{equation}\label{clco}
\co D^*f(\ox)(u)=\big\{A^*u\big|\;A\in\co\Bar\nabla f(\ox)\big\}\;\;\mbox{whenever}\;u\in\R^m.
\end{equation}

Now we recall the definition of a remarkable subclass of single-valued locally Lipschitzian mappings, which plays a very significant role in numerical optimization; see the books \cite{fp,is14} for the history and more discussions. This class can be described as follows. Given a mapping $f\colon\R^n\to\R^m$ locally Lipschitzian around $\ox$, we say that $f$ is {\em semismooth} at $\ox$ if
it is directionally differentiable at $\ox$ and the estimate
\begin{equation}\label{semi}
f(x)-f(\ox)-A(x-\ox)=o(\|x-\ox\|)
\end{equation}
holds when $x\to\ox$ and $A\in\co\Bar\nabla f(x)$. It is important to observe that estimate \eqref{semi} and the directional differentiability of $f$ are mutually independent assumptions; see, e.g., \cite{msz}.

The concept of semismoothness has been recently improved and extended in \cite{go19} to set-valued mappings. Recall from \cite{go19} that a set-valued mapping $F\colon\R^n\tto\R^m$ is {\em semismooth$^*$} at $(\ox,\oy)\in\gph F$ if whenever $(u,v)\in\R^n\times\R^m$ we have
\begin{equation}\label{semi6}
\la u^*,u\ra=\la v^*,v\ra\;\mbox{ for all }\;(v^*,u^*)\in\gph D^*F\big((\ox,\oy);(u,v)\big),
\end{equation}
where $D^*F((\ox,\oy);(u,v))$ is the {\em directional coderivative} of $F$ at $(\ox,\oy)\in\gph F$ in the direction $(u,v)\in\R^n\times\R^m$ defined by
\begin{equation}\label{dcod}
D^*F\big((\ox,\oy);(u,v)\big)(w):=\big\{q\in\R^n\big|\;(q,-w)\in N_{\scriptsize{\gph F}}\big((\ox,\oy);(u,v)\big)\big\}
\end{equation}
via the the {\em directional normal cone}
\begin{equation*}
N_{\scriptsize{\gph F}}\big((\ox,\oy);(u,v)\big)=\disp\Limsup_{\substack{t\dn 0\\(u',v')\to(u,v)}}\Hat N_{\scriptsize{\gph F}}\big((\ox,\oy)+t(u',v')\big)
\end{equation*}
to $\gph F$ at $(\ox,\oy)$ in the direction $(u,v)$. The latter construction was introduced in \cite{gin-mor}, while the directional coderivative \eqref{dcod} was defined and largely investigated in \cite{g}. As shown in \cite{go19}, for single-valued and locally Lipschitzian mapping $F=f$, the semismooth$^*$ property reduces to estimate \eqref{semi}, but without the directional differentiability requirement.\vspace*{0.03in}

The next proposition collects equivalent descriptions of semismoothness$^*$ for single-valued and locally Lipschitzian mappings that are utilized below.\vspace*{-0.07in}

\begin{Proposition}[\bf equivalent descriptions of semismoothness$^*$]\label{semi5} Let $f\colon\R^n\to\R^n$ be a locally Lipschitzian mapping around $\ox$. Then the following are equivalent:

{\bf(i)} $f$ is semismooth$^*$ at $\ox$.

{\bf(ii)} For any $x\to\ox$ and any $A\in\co\Bar\nabla f(x)$, we have estimate \eqref{semi}.

{\bf(iii)} For any $x\to\ox$ and any $w_x\in Df(x)(x-\ox)$ we have $f(x)-f(\ox)-w_x=o(\|x-\ox\|)$.\\
If in addition all the matrices $A\in\Bar\nabla f(x)$ are symmetric as $x$ sufficiently close to $\ox$, then the above conditions are equivalent to each of the listed below:

{\bf(iv)} For any $x\to\ox$ and any $w_x\in D^*f(x)(x-\ox)$ we have $f(x)-f(\ox)-w_x=o(\|x-\ox\|)$.

{\bf(v)} For any $x\to\ox$ and any $w_x\in\Hat D^*f(x)(x-\ox)$ we have $f(x)-f(\ox)-w_x=o(\|x-\ox\|)$.
\end{Proposition}\vspace*{-0.15in}
\begin{proof} The equivalence between (i) and (ii) is established in \cite[Proposition~3.7]{go19}. Since $Df(x)(x-\ox)\subset(\co\Bar\nabla f(\ox))(x-\ox)$, we get
(ii)$\implies$(iii). To justify (iii)$\implies$(ii), recall from \cite[Lemma~2.1]{s01} that (ii) is equivalent to the estimate
\begin{equation}\label{sun1}
f(x)-f(\ox)-\nabla f(x)(x-\ox)=o(\|x-\ox\|)\;\mbox{ for all }\;x\to\ox\;\mbox{ with }\;x\in\O_f.
\end{equation}
Suppose now that (iii) holds and pick $x\in\O_f$. Since $Df(x)(x-\ox)=\nabla f(x)(x-\ox)$, we get that (iii) implies \eqref{sun1}, and thus (ii) is satisfied.

Observe further that the additional symmetry assumption ensures that each $A\in\co\Bar\nabla f(x)$ is symmetric for $x$ near $\ox$, and thus (ii)$\implies$(iv) follows from \eqref{clco}. Implication (iv)$\implies$(v) is a direct consequence of the inclusion $\Hat D^*f(x)(x-\ox)\subset D^*f(x)(x-\ox)$. To verify finally (iv)$\implies$(ii), note that the symmetry of $\nabla f(x)$ for $x\in\O_f$ near $\ox$ tells us that $\Hat D^*f(x)(x-\ox)=\nabla f(x)(x-\ox)$ for such $x$ and hence (iv) yields \eqref{sun1}. Thus we deduce (ii) from \cite[Lemma~2.1]{s01}.
\end{proof}\vspace*{-0.07in}

It is worth mentioning that the symmetry assumption on $\Bar\nabla f(x)$ in Proposition~\ref{semi5} holds for important cases of mappings used in optimization. In particular, we have it for $f=\nabla\ph$, where $\ph$ is a function of class ${\cal C}^{1,1}$ near $\ox$. Indeed, it follows from \cite[Theorem~13.52]{rw} that $\Bar\nabla^2\ph(x):=\Bar\nabla(\nabla\ph)(x)$ is a compact set of symmetric matrices for such $x$. Furthermore, all the equivalences of Proposition~\ref{semi5} hold when $f:=P_r\ph$ is the proximal mapping of a prox-regular function $\ph$ with small $r>0$; see \eqref{proxmap}. This follows from \cite[Corollary~13.53]{rw}, which tells us that all the matrices in $\Bar\nabla(P_r\ph)(x)$ are symmetric in this case.\vspace*{0.03in}

For subsequent applications in this paper, we need to present some other notions and results of variational analysis. Recall that a mapping $f\colon\R^n\to\R^m$ is {\em semidifferentiable} at $\ox$ if there is a continuous and positively homogeneous operator $H\colon\R^n\to\R^m$ such that
\begin{equation*}
f(x)=f(\ox)+H(x-\ox)+o(\|x-\ox\|)\;\mbox{ for all }\;x\;\mbox{ near }\;\ox.
\end{equation*}
It follows from \cite[Exercise~9.25]{rw} that the semidifferentiability of $f$ implies that its graphical derivative $Df(\ox)$ is single-valued. By \cite[Proposition~9.50(b)]{rw} we have that semidifferentiability of locally Lipschitzian mappings is equivalent to their {\em proto-differentiability}, which means---in the general set-valued setting---that the graph of the mapping is derivable at the point in question. In the second-order framework of continuously prox-regular functions $\ph\colon\R^n\to\oR$ of our main interest here, the fundamental result of \cite[Theorem~13.40]{rw} tells us that the proto-differentiability of the subgradient mappings $\partial\ph$ is equivalent to the twice epi-differentiability of the function $\ph$. Recall to this end that the {\em second subderivative} of $\ph$ at $\ox\in\dom\ph$ for $\ov\in\R^n$ is the function $\d^2\ph(\ox,\ov)\colon\R^n\to[-\infty,\infty]$ defined for all $w\in\R^n$ by
\begin{equation*}
\d^2\ph(\bar x,\ov)(w):=\liminf_{t\dn 0,w'\to w}\Delta_t^2\ph(\bar x,\ov)(w')\;\mbox{ with }\;\Delta_t^2\ph(\bar x,\ov)(w):=\dfrac{\ph(\ox+tw')-\ph(\ox)-t\langle\ov,w'\rangle}{\frac{1}{2}t^2}.
\end{equation*}
According to \cite[Definition~7.23 and Proposition~7.2]{rw}, the {\em twice epi-differentiability} of $\ph$ at $\bar x$ for $\ov$ means that for every $w\in\R^n$ and every sequence $t_k\downarrow 0$ there exists a sequence $w_k\to w$ such that $\Delta_{t_k}^2\ph(\bar x,\ov)(w_k)\to\d^2\ph(\bar x,\ov)(w)$ as $k\to\infty$.\vspace*{0.03in}

To conclude this section, we show how to use the aforementioned results to derive a novel semismooth$^*$ inverse mapping theorem related to tilt-stable minimizers. As has been well recognized in optimization theory, inverse and implicit mapping theorems play a fundamental role in the justification of numerical algorithms; in particular, of the Newton type. Their proofs are usually involved while often using degree theory; see, e.g., \cite{fp} and the references therein. The proof of the following new result seems to be significantly easy even in standard settings.\vspace*{-0.07in}

\begin{Proposition}[\bf semismooth$^*$ inverse mapping theorem under prox-regularity and tilt stability]\label{eqiM} Let $\ph\colon\R^n\to\oR$ be continuously prox-regular at $\ox$ for $\ov=0$, where $\ox$ is a tilt-stable local minimizer of $\ph$ with some modulus $\kappa>0$. Then there are neighborhoods $U$ of $\ox$ and $V$ of $\ov$ such that the mapping $v\mapsto(\sub\ph)^{-1}(v)\cap U$ is single-valued and Lipschitz continuous on $V$. Furthermore, we have the equivalent assertions:

{\bf(i)} The mapping $\sub\ph$ is semismooth$^*$ at $(\ox,\ov)$, and $\ph$ twice epi-differentiable at $\ox$ for $\ov$.

{\bf(ii)} The mapping $v\mapsto(\sub\ph)^{-1}(v)\cap U$ is semismooth at $\ov$.
\end{Proposition}\vspace*{-0.15in}
\begin{proof} The existence of neighborhoods $U$ of $\ox$ and $V$ of $\ov$ for which the mapping $v\mapsto(\sub\ph)^{-1}(v)\cap U$ is single-valued and Lipschitz continuous on $V$
was established in \cite[Theorem~1.3]{pr98}. Denote $g(v):=(\sub\ph)^{-1}(v)\cap U$ on $V$ and observe that the mapping $\sub\ph$ is semismooth$^*$ at $(\ox,\ov)$ if and only if $g$ is semismooth$^*$ at $\ov$. It follows from the continuous prox-regularity of $\ph$ at $\ox$ for $\ov$ due to the aforementioned result of \cite[Theorem~13.40]{rw} that the twice epi-differentiability of $\ph$ at $\ox$ for $\ov$ amounts to the proto-differentiability of $\sub\ph$ therein. This is equivalent to the proto-differentiability of $g$ at $\ov$ for $\ox$ and also, due to the local Lipschitz continuity of $g$, to the semidifferentiability of $g$ at $\ov$. Furthermore, the result of \cite[Proposition~2D.1]{dr} tells us that the latter property reduces in this setting to the classical directional differentiability of $g$ at $\ov$. Employing finally Proposition~\ref{semi5} verifies the claimed equivalence between assertions (i) and (ii) of the theorem, where the semismooth and semismooth$^*$ properties of $g$ are the same due to the established directional differentiability of $g$ at $\ov$ in this case. \end{proof}\vspace*{-0.25in}

\section{Basic Estimates for Newton Iterations}\sce\label{sect03}\vspace*{-0.05in}

In this section we derive some technical estimates, which play a crucial role in the subsequent justification of both coderivative-based and graphical derivative-based  generalized Newton algorithms. It is important to emphasize that our results provide not only {\em qualitative} but also {\em quantitative} estimates involving moduli of tilt stability.\vspace*{0.03in}

We start with the coderivative estimates for tilt-stable local minimizers of ${\cal C}^{1,1}$ functions.\vspace*{-0.07in}

\begin{Theorem}[\bf estimates of Newton iterations via coderivatives]\label{codne} Let $\ph\colon\R^n\to\oR$ be a ${\cal C}^{1,1}$ function on a neighborhood of $\ox$, and let $\ox$ be a tilt-stable local minimizer for $\ph$ with modulus $\kappa>0$. Then there exists $\dd>0$ such that for any $(x,v)\in(\gph\nabla\ph)\cap\B_\dd(\ox,\ov)$ with $\ov=0$ and any $q\in(D^*\nabla\ph)^{-1}(v,x)(v-\ov)$ we find $y_x\in(D^*\nabla\ph)(x)(x-\ox)$ satisfying
\begin{equation}\label{inv90}
\|x-\ox-q\|\le\kappa\,\|\nabla\ph(x)-\nabla\ph(\ox)-{y}_x\|.
\end{equation}
\end{Theorem}\vspace*{-0.15in}
\begin{proof} Let us first derive a similar estimate for vectors $q$ belonging to the regular coderivative \eqref{rcod} of $\nabla\ph^{-1}$. Namely, take any $q\in(\Hat D^*\nabla\ph)^{-1}(v,x)(v-\ov)$ in the setting of the theorem and show that there exists $\tilde{y}_x\in(D^*\nabla\ph)(x)(x-\ox)$ such that
\begin{equation}\label{inv03}
\|x-\ox-q\|\le\kappa\,\|w-\nabla\ph(\ox)-\tilde{y}_x\|.
\end{equation}
To verify \eqref{inv03}, recall the uniform second-order growth characterization \eqref{usogc} of the tilt-stable local minimizer $\ox$ obtained in Proposition~\ref{usogc2}. It gives us neighborhoods $U$ of $\ox$ and $V$ of $\ov$ with
\begin{equation}\label{ugc}
\la x-x',v-v'\ra\ge\kappa^{-1}\|x-x'\|^2\;\mbox{ whenever }\;(x,v),\,(x',v')\in(\gph\nabla\ph)\cap(U\times V).
\end{equation}
Suppose with no harm that $\nabla\ph$ is Lipschitz continuous on $U$ with constant $\ell>0$ and that $U\subset\B_\eta(\ox)$, where $\eta$ is taken from estimate \eqref{lip} in Proposition~\ref{lipcon}
with $f:=\nabla\ph$. Let $\dd>0$ be so small that $\B_\dd(\ox,\ov)\subset{U}\times{V}$ and then pick $(x,v)\in(\gph\nabla\ph)\cap\B_\dd(\ox,\ov)$ and $q\in(\Hat D^*\nabla\ph)^{-1}(v,x)(w-\ov)$ with $w\in\R^n$. We get $\ov-w\in(\Hat D^*\nabla\ph)(x)(-q)$, which tells us by \eqref{rcod} that $(\ov-w, q)\in\Hat N_{\scriptsize{\gph\nabla\ph}}(x,v)$. Using \eqref{rnc} implies that for any $\epsilon>0$ there is $r>0$ with
\begin{equation}\label{sm}
\la\ov-w,u-x\ra+\la q,p-v\ra\le\epsilon\big(\|u-x\|+\|p-v\|\big)\;\mbox{ when }\;(u,p)\in(\gph\nabla\ph)\cap\B_{r}(x,v).
\end{equation}
Suppose further without loss of generality that $\B_r(x,v)\subset U\times V$. Define $u_{t}:=x+t(q-x+\ox)$ and $p_{t}:=\nabla\ph(u_{t})$ for $t>0$ and then observe that $(u_{t},p_{t})\in\gph\nabla\ph$ and that $(u_{t},p_{t})\to(x,v)$ as $t\downarrow 0$. This allows us to obtain the inclusions
\begin{equation*}
(u_{t},p_{t})\in(\gph\nabla\ph)\cap\B_{r}(x,v)\subset(\gph\nabla\ph)\cap(U\times V)\;\mbox{ for all }\;t>0.
\end{equation*}
Employing now \eqref{ugc} brings us to the conditions
\begin{equation*}
\la u_{t}-x,p_{t}-v\ra\ge\kappa^{-1}\|u_{t}-x\|^2=\kappa^{-1}t^2\|q-x+\ox\|^2,
\end{equation*}
which together with $q=t^{-1}(u_{t}-x)+x-\ox$ provide the estimate
\begin{equation}\label{sm01}
\begin{array}{lll}
\la q,p_{t}-v\ra&=&t^{-1}\la u_{t}-x,p_{t}-v\ra+\la x-\ox,p_{t}-v\ra\\
&\ge&\kappa^{-1}t\|q-x+\ox\|^2+\la x-\ox,p_{t}-v\ra.
\end{array}
\end{equation}
Plugging $(u_{t},p_{t})$ into \eqref{sm} and appealing to \eqref{sm01} lead us to
\begin{eqnarray}\label{sm001}
\kappa^{-1}t\|q-x +\ox\|^2&\le& \la q,p_{t}-v\ra+\la x-\ox,v-p_{t}\ra\nonumber\\
&\le&\epsilon\big(t\|q-x+\ox\|+\|\nabla\ph(u_t)-\nabla\ph(x)\|\big)\nonumber\\
&&+\la w-\ov,x+t(q_v-x+\ox)-x\ra+\la x-\ox,v-p_{t}\ra\nonumber\\
&\le&t\epsilon\big(\|q-x+\ox\|+\ell\|q-x+\ox\|\big)\nonumber\\
&&-t\la w-\ov,x-\ox-q\ra+\la x-\ox,\nabla\ph(x)-\nabla\ph(u_t)\ra,
\end{eqnarray}
where $\ell$ is a Lipschitz constant of $\nabla\ph$ on $U$. It follows from $[u_{t},x]\subset U$ that the function $z\mapsto\la x-\ox,\nabla\ph\ra(z)$ is Lipschitz continuous on an open set containing $[u_{t},x]$. Applying now the mean value inequality from \cite[Corollary~4.14(ii)]{m18} to the latter function and using the coderivative scalarization formula from \cite[Theorem~1.32]{m18} give us  vectors $c_t\in[u_{t},x)$ and $y_t\in\sub\la x-\ox,\nabla\ph\ra(c_t)=(D^*\nabla\ph)(c_t)(x-\ox)$ that satisfy the conditions
\begin{equation*}
\la x-\ox,\nabla\ph(x)-\nabla\ph(u_{t})\ra\le\la y_t,x-u_t\ra=t\la y_t,x-\ox-q\ra.
\end{equation*}
Using them together with (\ref{sm001}), we arrive at the estimate
\begin{equation}\label{sm002}
\|q-x+\ox\|^2\le\kappa\epsilon(1+\ell)\|q-x+\ox\|+\kappa\la y_t-w+\ov,x-\ox-q\ra.
\end{equation}
Since $y_t\in(D^*\nabla\ph)(c_t)(x-\ox)$ with $c_t\in U$, it follows from the above choice of $U\subset\B_\eta(\ox)$ and Proposition~\ref{lipcon} that there exists a positive number $\rho$ such that  $\|y_t\|\le\rho\,\|x-\ox\|$. This allows us to claim without loss of generality the existence of $\tilde{y}_x\in\R^n$ such that $y_t\to\tilde{y}_x$ as $t\dn 0$. Observing that $c_t\to x$ as $t\dn0$, we get $\tilde{y}_x\in(D^*\nabla\ph)(x)(x-\ox)$. Furthermore, the passage to the limit in \eqref{sm002} as $t\dn 0$ gives us the estimate
\begin{equation*}
\|q-x+\ox\|^2\le\kappa\epsilon(1+\ell)\|q-x+\ox\|+\kappa\la\tilde{y}_x-w+\ov,x-\ox-q\ra.
\end{equation*}
Finally, letting $\epsilon\dn 0$ brings us to
\begin{eqnarray*}
\|x-\ox-q\|^2\le\kappa\la\tilde{y}_x-w+\ov,x-\ox-q\ra&\le&\kappa\|w-\ov-\tilde{y}_x\|\cdot\|x-\ox-q\|\\
&=&\kappa\|w-\nabla f(\ox)-\tilde{y}_x\|\cdot \|x-\ox-q\|,
\end{eqnarray*}
which proves \eqref{inv03} for the vector $q\in(\Hat D^*\nabla\ph)^{-1}(v,x)(v-\ov)$ from the regular coderivative.\vspace*{0.03in}

Our next step is to verify estimate \eqref{inv90} for any selected vector $q\in(D^*\nabla\ph)^{-1}(v,x)(v-\ov)$ from the limiting coderivative \eqref{lcod}. To proceed, take the number $\dd>0$ for which we derived \eqref{inv03} and then pick any $(x,v)\in(\gph\nabla\ph)\cap\B_{\dd/2}(\ox,\ov)$ and $q\in(D^*\nabla\ph)^{-1}(v,x)(v-\ov)$. This implies by the coderivative and normal cone definitions the existence of sequences $(v_k,x_k)\to(v,x)$ and $(w_k,q_k)\to(v-\ov,q)$ as $k\to\infty$ with $q_k\in(\Hat D^*\nabla\ph)^{-1}(v_k,x_k)(w_k)$. Remember that $\ov=0$ and thus $q_k\in(\Hat D^*\nabla\ph)^{-1}(v_k,x_k)(w_k-\ov)$. Employing \eqref{inv03} tells us that for all $k$ sufficiently large there exists $y_k\in(D^*\nabla\ph)(x_k)(x_k-\ox)$ such that the estimate
\begin{equation}\label{cod-est}
\|x_k-\ox-q_k\|\le\kappa\,\|w_k-\nabla\ph(\ox)-y_k\|
\end{equation}
holds. By Proposition~\ref{lipcon} and the fact that $x_k\to x$, we can assume with no harm that the sequence $\{y_k\}$ is bounded, and hence there exists its subsequence that converges to some $y_x\in\R^n$. Passing finally to the limit in \eqref{cod-est} as $k\to\infty$ verifies \eqref{inv90} and thus completes the proof.
\end{proof}\vspace*{-0.07in}

Next we intend to derive certain counterparts of Theorem~\ref{codne} with employing the graphical derivative. Due the absence (to the best of our knowledge) an appropriate mean value theorem for the graphical derivative, we need to either impose an additional semidifferentiability assumption, or to use the extended Hessian set
\begin{equation}\label{hes}
\co\Bar\nabla^2\ph(x):=\co\Bar\nabla(\nabla\ph)(x),\quad x\in\R^n.
\end{equation}\vspace*{-0.35in}

\begin{Theorem}[\bf Newton estimates involving graphical derivatives]\label{inv02} Let $\ph\colon\R^n\to\oR$ be a ${\cal C}^{1,1}$ function on a neighborhood of $\ox$, which is a tilt-stable local minimizer for $\ph$ with modulus $\kappa>0$, and let $\ov=0$. Then the following assertions hold:

{\bf(i)} There exists $\dd>0$ such that for any $(x,v)\in(\gph\nabla\ph)\cap\B_\dd(\ox,\ov)$ and $q\in(D\nabla\ph)^{-1}(v,x)(v-\ov)$ we find a matrix $A\in\co\Bar\nabla^2\ph(x)$ for which the estimate
\begin{equation}\label{gdi}
\|x-\ox-q\|\le\kappa\,\|\nabla\ph(x)-\nabla\ph(\ox)-A(x-\ox)\|
\end{equation}
is satisfied. Similarly, there exists $\dd>0$ such that for any $(x,v)\in(\gph\nabla\ph)\cap\B_\dd(\ox,\ov)$ and $q\in(D\nabla\ph)^{-1}(v,x)(\ov-v)$ we find a matrix $A\in\co\Bar\nabla^2\ph(x)$ ensuring the estimate
\begin{equation}\label{gdi3}
\|x-\ox+q\|\le\kappa\,\|\nabla\ph(x)-\nabla\ph(\ox)-A(x-\ox)\|.
\end{equation}

{\bf(ii)} If in addition the gradient mapping $x\mapsto\nabla\ph(x)$ is semidifferentiable on a neighborhood of $\ox$, then there exists $\dd>0$ such that for any $(x,v)\in(\gph\nabla\ph)\cap\B_\dd(\ox,\ov)$
and any $q\in(D\nabla\ph^{-1})(v,x)(v-\ov)$ we have the estimate
\begin{equation}\label{gdi1}
\|x -\ox-q\|\le\kappa\,\|\nabla\ph(x)-\nabla\ph(\ox)-(D\nabla\ph)(x)(x-\ox)\|.
\end{equation}
\end{Theorem}\vspace*{-0.15in}
\begin{proof} Let us begin with verifying (ii). The assumed tilt stability of $\ox$ gives us neighborhoods $U$ of $\ox$ and $V$ of $\ov$ for which \eqref{ugc} holds. We also suppose that $U$ is the neighborhood on which the mapping $x\mapsto\nabla\ph(x)$ is semidifferentiable. Take $\dd>0$ such that $\B_\dd(\ox,\ov)\subset U\times V$ and then pick $(x,v)\in(\gph\nabla\ph)\cap\B_\dd(\ox,\ov)$ and $q\in(D\nabla\ph)^{-1}(v,x)(v-\ov)$. The latter implies that $v-\ov\in(D\nabla\ph)(x)(q)$, and so $(q,v-\ov)\in T_{\scriptsize{\gph\nabla\ph}}(x,v)$ by \eqref{gder}. It follows from the tangent cone definition \eqref{tan} that there exist sequences $t_k\dn 0$ and $(q_k,w_k)\to(q,v-\ov)$ as $k\to\infty$ with $(x+t_k q_k,v+t_k w_k)\in\gph\nabla\ph$, $k\in\N$. Denoting $u_k:=x+t_k(x-\ox)$ and $z_k:=\nabla\ph(u_k)$, we get get for all large $k$ that $(u_k,z_k)\in(\gph\nabla\ph)\cap(U\times V)$. Then \eqref{ugc} tells us that
\begin{equation*}
\la u_k-(x+t_k q_k),z_k-(v+t_k w_k)\ra\ge\kappa^{-1}\|u_k-(x+t_k q_k)\|^2
\end{equation*}
when the index $k$ is sufficiently large. This ensures the estimate
\begin{equation}\label{fj01}
t_k\|x-\ox-q_k\|^2\le\kappa\,\big\la x-\ox-q_k,\nabla\ph(u_k)-\nabla\ph(x)-t_k w_k\big\ra.
\end{equation}
Using the semidifferentiability of $\nabla\ph$ at any $x\in U$, we deduce from \cite[Exercise~9.25]{rw} that the mapping $(D\nabla\ph)(x)$ is single-valued with
\begin{equation*}
\nabla\ph\big(x+t_k(x-\ox)\big)=\nabla\ph(u_k)=\nabla\ph(x)+t_k(D\nabla\ph)(x)(x-\ox)+o(t_k)
\end{equation*}
for all large $k$. Combining it with \eqref{fj01} yields the inequality
\begin{equation*}
\|x-\ox-q_k\|^2\le\kappa\|x-\ox-q_k\|\cdot\Big\|(D\nabla\ph)(x)(x-\ox)+\frac{o(t_k)}{t_k}-w_k\Big\|.
\end{equation*}
Passing there to the limit as $k\to\infty$ gives us \eqref{gdi1} and completes the proof of (ii).\vspace*{0.03in}

Next we justify assertion (i) while focusing on estimate \eqref{gdi} with $q\in(D\nabla\ph^{-1})(v,x)(v-\ov)$. To verify first the one in \eqref{fj01}, employ the mean value theorem from \cite[Proposition~7.1.16]{fp} together with the classical Carath\'eodory theorem and find in this way $y_i^k\in (u_k,x)$ and $\al_i^k\ge 0$ as $i=1,\ldots,n^2+1$ with $\sum_{i=1}^{n^2+1}\al_i^k=1$ so that
\begin{equation*}
\nabla\ph(u_k)-\nabla\ph(x)=t_k\sum_{i=1}^{n^2+1}\al_i^k A_i^k(x-\ox)\;\mbox{ with }\;A_i^k\in\co\big(\Bar\nabla^2\ph(y_i^k)\big).
\end{equation*}
Letting $k\to\infty$ tells us that $y_i^k\to x$ for any $i=1,\ldots,n^2+1$. Since $\ox$ is a tilt-stable local minimizer of $\ph$, it follows from \cite[Theorem~2.1]{pr98} that
\begin{equation*}
\la w,u\ra>0\;\mbox{ for all }\;w\in(D^*\nabla\ph)(\ox)(u),\;u\ne 0.
\end{equation*}
By \eqref{clco} the latter amounts to saying that
\begin{equation}\label{nmc2}
\la Au,u\ra>0\;\mbox{ whenever }\;0\ne u\in\R^n,\;A\in\co\Bar\nabla^2\ph(\ox).
\end{equation}
This implies that all the matrices from the Hessian set $\co\Bar\nabla^2\ph(\ox)$ are nonsingular. Employing \cite[Lemma~7.5.2]{fp}, suppose without loss of generality that the sequence $\{A_k\}$ is bounded. Passing to a subsequence if necessary tells us that $\al_i^k\to\al_i$ for some $\al_i\ge 0$ with $\sum_{i=1}^{n^2+1}\al_i=1$ and that $A_i^k\to A_i$ for some $n\times n$ matrix $A_i$ as $k\to\infty$ whenever $i=1,\ldots,n^2+1$. Then \cite[Proposition~7.1.4]{fp} yields the existence of $A_i\in\co(\Bar\nabla^2\ph(x))$. Define $A:=\sum_{i=1}^{n^2+1}\al_i A_i$ and observe that $A\in\co(\Bar\nabla^2\ph(x))$. Combining it with \eqref{fj01} and passing to the limit as $k\to\infty$ bring us to
\begin{equation*}
\|x-\ox-q\|^2\le\kappa\la x-\ox-q,A(x-\ox)-v+\ov\ra\le\kappa\|x-\ox-q\|\cdot\|v-\ov-A(x-\ox)\|,
\end{equation*}
which justifies the claimed estimate \eqref{gdi}. Finally, estimate \eqref{gdi3} can be justified similarly by choosing $u_k$ in the proof of \eqref{gdi} as $u_k:=x-t_k(x-\ox)$
and then proceeding as above.
\end{proof}\vspace*{-0.25in}

\section{Coderivative-Based Newton Algorithm for ${\cal C}^{1,1}$ Functions}\sce\label{sect04}\vspace*{-0.05in}

This section is devoted to the design and justification of a new coderivative-based generalized Newton algorithm for title-stable minimizers of ${\cal C}^{1,1}$ functions. Given such a function $\ph\colon\R^n\to\R$ around its tilt-stable local minimizer $\ox$ in the unconstrained problem \eqref{main}, define the set-values mapping $\Upsilon^*\colon\R^n\rightrightarrows\R^n$ by
\begin{equation}\label{di1}
\Upsilon^*(x):=\big(D^*\nabla\ph\big)^{-1}\big(\nabla\ph(x),x\big)\big(\nabla\ph(x)\big)=\big\{y\in\R^n\big|-\nabla\ph(x)\in(D^*\nabla\ph)(x)(-y)\big\}.
\end{equation}
Note that the coderivative in \eqref{di1} admits the aforementioned subdifferential scalarization
\begin{equation}\label{scal}
D^*\nabla\ph(x)(-y)=\partial\la-y,\nabla\ph\ra(x)\;\mbox{ for }\;x\;\mbox{ near }\;\ox,
\end{equation}
while the convex hull of the set on the right-hand side agrees with Clarke's generalized gradient; see \cite{m18,rw}. Representation \eqref{scal} significantly simplifies the computation in \eqref{di1}. When $\ph$ is ${\cal C}^2$-smooth around $\ox$, the set $\Upsilon^*(x)$ in \eqref{di1} reduces to $\nabla^2\ph(x)^{-1}\nabla\ph(x)$ for all $x$ near $\ox$ while resembling the directions in the classical Newton method.\vspace*{0.03in}

Using \eqref{di1}, we now propose the following Newton-type algorithm for ${\cal C}^{1,1}$ functions.\vspace*{-0.07in}

\begin{Algorithm}[\bf generalized Newton method for ${\cal C}^{1,1}$ functions via coderivatives]\label{nm1} {\rm Pick $x_0\in\R^n$ and set $k:=0$.\\
{\bf Step~1:} If $\nabla\ph(x_k)=0$, then stop.\\
{\bf Step~2:} Otherwise, select a direction $d_k\in\Upsilon^*(x_k)$ and set $x_{k+1}:=x_k-d_k$.\\
{\bf Step~3:} Let $k\leftarrow k+1$ and then go to Step~1.}
\end{Algorithm}\vspace*{-0.05in}

To proceed with the study of Algorithm~\ref{nm1}, first we should clarify the {\em solvability} issue. It is done in the next proposition for the case of tilt-stable minimizers of our main interest.\vspace*{-0.07in}

\begin{Proposition}[\bf solvability of subproblems  the coderivative-based Newton algorithm]\label{exist} Let $\ph\colon\R^n\to\R$ be a ${\cal C}^{1,1}$ function on a neighborhood of $\ox$, which is a tilt-stable local minimizer for $f$ with some modulus $\kappa>0$. Then there exists a neighborhood $O$ of $\ox$ such that the set-valued mapping $\Upsilon^*(x)$ from \eqref{di1} is nonempty and compact-valued for all $x$ in $O$.
\end{Proposition}\vspace*{-0.2in}
\begin{proof} Since $\ox$ is a tilt-stable minimizer of $\ph$, we conclude from \cite[Theorem~1.3]{pr98} that there are neighborhoods $U$ of $\ox$ and $V$ of $\ov=0$ such that
the mapping $v\mapsto(\nabla\ph)^{-1}(v)\cap U$ is Lipschitz continuous on $V$. By the Lipschitz continuity of $\nabla\ph$ around $\ox$ we find numbers $\dd_1,\dd_2>0$ with
\begin{equation*}
\B_{\dd_1}(\ox)\times\B_{\dd_2}(\ov)\subset{U}\times{V}\;\mbox{ and }\;\nabla\ph(x)\in\B_{\dd_2}(\ov)\;\mbox{ for all }\;x\in\B_{\dd_1}(\ox).
\end{equation*}
Pick further $x\in\B_{\dd_1}(\ox)$ and denote $v:=\nabla\ph(x)$, which implies that $(x,v)\in U\times V$. Considering the mapping $g(v):=(\nabla\ph)^{-1}(v)\cap U$ on $V$, observe that it is locally Lipschitzian around $v$. It follows from the above constructions and the scalarization formula that
\begin{equation*}
(D^*\nabla\ph)^{-1}(v,x)(u)=D^*g(v)(u)=\partial\la u,g\ra(v)\;\mbox{ for any }\;u\in\R^n.
\end{equation*}
and so the set $\Upsilon^*(x)$ is nonempty and compact for all $x$ near $\ox$ due to, e.g., \cite[Theorem~1.22]{m18}.
\end{proof}\vspace*{-0.07in}

Note that the mapping $\Upsilon^*$ is neither a Newton approximation in the sense of \cite{fp}, nor a Newton map in the sense of \cite{kk}. The latter map is in fact a collection of linear mappings.\vspace*{0.03in}

We are now in a position to establish superlinear convergence of the Newton method from Algorithm~\ref{nm1}. Recall that a sequence $\{x_k\}$ converging to $\ox$ converges {\em Q-superlinearly} if $\|x_{k+1}-\ox\|=o(\|x_k-\ox\|)$ as $k\to\infty$. In what follows we drop the letter $Q$ and simply speak about superlinear convergence of a sequence.\vspace*{-0.07in}

\begin{Theorem}[\bf superlinear convergence of the coderivative-based Newton algorithm for ${\cal C}^{1,1}$ functions]\label{des12} Let $\ph\colon\R^n\to\R$ be a ${\cal C}^{1,1}$ function on a neighborhood of its tilt-stable local minimizer $\ox$ with modulus $\kappa>0$, and let the gradient mapping $\nabla\ph$ be semismooth$^*$ at $\ox$. Then there is $\dd>0$ such that for any starting point $x_0\in\B_\dd(\ox)$ we get that every sequence $\{x_k\}$ constructed by Algorithm~{\rm\ref{nm1}} converges to $\ox$ and the rate of convergence is superlinear.
\end{Theorem}\vspace*{-0.15in}
\begin{proof} Since $\ph$ is a ${\cal C}^{1,1}$ function on a neighborhood of $\ox$, it follows from \cite[Theorem~13.52]{rw} that all the matrices in $\Bar\nabla^2\ph(x)$ are symmetric for any $x$ close to $\ox$. The imposed semismooth$^*$ assumption on $\nabla\ph$ at $\ox$ gives us by Proposition~\ref{semi5}(iv) a positive number $\dd'$ such that
\begin{equation}\label{mss2}
\|\nabla\ph(x)-\nabla\ph(\ox)-y_x\|<\frac{1}{2\kappa}\|x-\ox\|\;\mbox{ for all }\;x\in\B_{\dd'}(\ox)\;\mbox{ and }\;\;y_x\in(D^*\nabla\ph)(x)(x-\ox).
\end{equation}
Take further $\dd>0$ from Theorem~\ref{codne} for which estimate \eqref{inv90} holds. Proposition~\ref{exist} ensures the existence of a neighborhood $O$ of $\ox$ on which $\Upsilon^*(x)\ne\emp$. Since $\nabla\ph$ is Lipschitz continuous around $\ox$, we find positive numbers $\dd_1<\dd'$ and $\dd_2$ such that $\B_{\dd_1}(\ox)\times\B_{\dd_2}(\ov)\subset\B_\dd(\ox,\ov)$, $\B_{\dd_1}(\ox)\subset O$,
and $\nabla\ph(x)\in\B_{\dd_2}(\ov)$ for all $x\in\B_{\dd_1}(\ox)$.

Pick now $x_k\in\B_{\dd_1}(\ox)$ for $k$ sufficiently large and suppose without loss of generality that it holds for all $k\in\N$. Then we get that $\Upsilon^*(x_k)\ne\emp$. According to Algorithm~\ref{nm1}, select a direction $d_k\in\Upsilon^*(x_k)$ and set $x_{k+1}:=x_k-d_k$. Denoting $v_k:=\nabla\ph(x_k)$ ensures that $(x_k,v_k)\in(\gph\nabla\ph)\cap\B_\dd(\ox,\ov)$. Since $\ov=0$, we have $d_k\in(D^*\nabla\ph)^{-1}(v_k,x_k)(v_k-\ov)$ and then deduce from Theorem~\ref{codne} that there exists a vector $y_k\in(D^*\nabla\ph)(x_k)(x_k-\ox)$ providing the estimate
\begin{equation}\label{newt2}
\|x_k-\ox-d_k\|\le\kappa\,\|\nabla\ph(x_k)-\nabla\ph(\ox)-y_k\|.
\end{equation}
Combining the latter with $x_k\in\B_{\dd_1}(\ox)\subset\B_{\dd'}(\ox)$ and \eqref{mss2}, we obtain that
\begin{equation*}
\|x_{k+1}-\ox\|<\sm\|x_k-\ox\|,
\end{equation*}
which yields $x_{k+1}\in\B_{\dd_1}(\ox)$. This tells us that every sequence $\{x_k\}$ with the starting point $x_0\in\B_{\dd_1}(\ox)$ generated by Algorithm~\ref{nm1} is contained in $\B_{\dd_1}(\ox)$ and converges to $\ox$ as $k\to\infty$. Employing this, Proposition~\ref{semi5}(iv), and estimate \eqref{newt2} ensures that
\begin{equation*}
\|x_{k+1}-\ox\|=o(\|x_k-\ox\|)
\end{equation*}
for all large $k$, and hence it shows that the rate of convergence of $x_k\to\ox$ is superlinear.
\end{proof}\vspace*{-0.17in}

\begin{Remark}[\bf comparison with the semismooth Newton method]\label{semi-comp} {\rm To compare Algorithm~\ref{nm1} with the celebrated semismooth method to solve problems \eqref{main} of unconstrained optimization with ${\cal C}^{1,1}$ objectives, consider the mapping $\Bar\Upsilon\colon\R^n\rightrightarrows\R^n$ defined by
\begin{equation*}
\Bar\Upsilon_{\tiny \mbox{C}}(x):=\big\{d\in\R^n\big|\;\nabla\ph(x)=Ad\;\;\mbox{with some}\;\;A\in\co\Bar\nabla^2\ph(x)\big\},
\end{equation*}
where the set $\co\Bar\nabla^2\ph(x)$ is taken from \eqref{hes}. Then the semismooth Newton algorithm for \eqref{main}, induced by the usage of Clarke's generalized Jacobian to solve semismooth equations, is formulated similarly to Algorithm~\ref{nm1} with the replacement of $\Upsilon^*(x)$ therein by $\Bar\Upsilon_{\tiny \mbox{C}}(x)$; see, e.g., \cite[Algorithm~2.56]{is14} and \cite[Algorithm~7.5.1]{fp} for more details. As an advantage of Algorithm~\ref{nm1} over the semismooth method, we mention a better coderivative calculus and the subdifferential representation \eqref{scal}, which is is not available for Clarke's constructions. Note also that for ${\cal C}^{1,1}$ functions we get from Proposition~\ref{semi5}(ii) that the semismooth$^*$ property and estimate \eqref{semi} in the standard semismooth property are equivalent. Furthermore, the best known superlinear convergence result presented, e.g., in \cite[Theorem~2.57]{is14} establishes superlinear convergence of the semismooth
Newton method under the second-order condition \eqref{nmc2}. As follows from \eqref{clco}, the latter condition is equivalent to tilt stability of the local minimizer $\ox$ utilized in Theorem~\ref{nm1}.
The only difference between these two algorithms is the set of eligible directions at each iteration. Indeed, \eqref{clco} tells us that the inclusion $\Upsilon^*(x)\subset\Bar\Upsilon_{\tiny \mbox{C}}(x)$ holds for all $x\in \R^n$.

Observe that it is possible to replace the second-order condition \eqref{nmc2} by
\begin{equation*}
\la Hu, u\ra>0\quad\mbox{whenever}\;\;0\ne u\in\R^n,\;H\in\Bar\nabla^2\ph(\ox).
\end{equation*}
This condition is clearly equivalent to \eqref{nmc2}. Defining the mapping $\Bar\Upsilon\colon\R^n\rightrightarrows\R^n$ by
\begin{equation}\label{dib}
\Bar\Upsilon(x):=\big\{d\in\R^n\;\big|\;\nabla\ph(x)=H d\;\;\mbox{for some}\;\;H\in\Bar\nabla^2\ph(x)\big\},\;\;x\in\R^n,
\end{equation}
we can get a Newton method for the optimization problem \eqref{main} by replacing $\Upsilon^*$ in Algorithm~\ref{nm1} by $\Bar\Upsilon$. This was championed by Qi  in \cite{qi93} for equations. The main disadvantage of latter method is the difficulty of calculating of $\Bar\nabla^2\ph$. Since for any $x\in\R^n$ the inclusion $\Bar\nabla^2\ph(x)(d)\subset(D^*\nabla\ph)(x)(d)$ always holds, we readily arrive at the inclusion $\Bar\Upsilon(x)\subset\Upsilon^*(x)$.}
\end{Remark}

Next we present an example showing that the semismooth$^*$ assumption in Theorem~\ref{des12} is essential for the convergence of Algorithm~\ref{nm1} under the fulfillment of all the other assumptions of the theorem. This example is based on \cite[Example~2.4]{jqcs} for the case of equations.\vspace*{-0.05in}

\begin{Example}[\bf failure of convergence of Newton iterations in the absence of semismooth$^*$ property]\label{exa} {\rm Define the Lipschitz continuous function $\psi\colon\R\to\R$ by
\begin{equation*}
\psi(x):=\left\{\begin{array}{ll}
x^2\sin\frac{1}{x}+2x&{\rm if}\;\;x\ne 0,\\
0&{\rm if}\;\;x=0.
\end{array}
\right.
\end{equation*}
It it shown in \cite{jqcs} for the semismooth Newton method that estimate \eqref{semi} fails for the function $\psi$ at $\ox=0$, and that the semismooth Newton iterations for solving the equation $\psi(x)=0$ starting with $x_0=\frac{1}{2\alpha\pi}$, for any fixed number $\alpha>0$, do not converge to $\ox$. Proposition~\ref{semi5} tells us that the semismooth$^*$ property of $\psi$ at $\ox$ fails as well. Consider now the function
\begin{equation*}
\ph(x):=\int_{0}^{x}\psi(t)dt,\quad x\in\R,
\end{equation*}
for which $\nabla\ph(x)=\psi(x)$ on $\R$. It is easy to check that
\begin{equation*}
\co\Bar\nabla^2\ph(\ox)=\co\Bar\nabla\psi(\ox)=[1,3]\;\mbox{ with }\;\nabla\ph(\ox)=\psi(\ox)=0,
\end{equation*}
and that condition \eqref{nmc2} is satisfied. This tells us that $\ox$ is a tilt-stable local minimizer of $\ph$. Since $\psi$ is continuously differentiable at every point but the origin, all the generalized derivatives for $\psi$ at $x\ne 0$ reduce to the classical one. Combining all of this with the result of \cite[Example~2.4]{jqcs} shows that the iterations of Algorithm~\ref{nm1} starting with $x_0=\frac{1}{2\alpha\pi}$ do not converge to the tilt-stable local minimizer in question.}
\end{Example}
\vspace*{-0.2in}

\section{${\cal C}^{1,1}$-Newton Algorithm Based on Graphical Derivatives}\sce\label{sect04a}\vspace*{-0.05in}

In this section we develop another Newton-type algorithm for tilt-stable local minimizers of ${\cal C}^{1,1}$ functions in \eqref{main}. The difference between Algorithm~\ref{nm1} and the new one is that now we are based on graphical derivatives \eqref{gder} of gradient mappings instead of coderivatives as in Section~\ref{sect04}. As shown in this section, assuming further that the objective function $\ph$ is twice epi-differentiable allows us to design a new algorithm with subproblems whose optimal solutions can be used to update the proposed algorithm.\vspace*{0.03in}

Let $\ph\colon\R^n\to\R$ be a ${\cal C}^{1,1}$ function around a point $\ox$. For all $x$ near $\ox$ define the sets
\begin{equation}\label{di1g}
\Upsilon(x):=\big(D\nabla\ph\big)^{-1}\big(\nabla\ph(x),x\big)\big(-\nabla\ph(x)\big)=\big\{y\in\R^n\big|\;-\nabla\ph(x)\in\big(D\nabla\ph\big)(x)(y)\big\}.
\end{equation}
Using these sets, we formulate now the following Newton-type algorithm.\vspace*{-0.05in}

\begin{Algorithm}[\bf generalized Newton method for ${\cal C}^{1,1}$ functions via graphical derivatives]\label{nm2} {\rm Pick $x_0\in\R^n$ and set $k:=0$.\\
{\bf Step~1:} If $\nabla\ph(x_k)=0$, then stop.\\
{\bf Step~2:} Otherwise, select a direction $d_k\in\Upsilon(x_k)$ and set $x_{k+1}:=x_k+d_k$.\\
{\bf Step~3:} Let $k\leftarrow k+1$ and then go to Step~1.}
\end{Algorithm}\vspace*{-0.05in}

To proceed further, first we have to address the solvability of subproblems  in Algorithm~\ref{nm2}, which is resolved in the next proposition provided that the function $\ph$ is of class ${\cal C}^{1,1}$ in a neighborhood of its tilt-stable local minimizer.\vspace*{-0.07in}

\begin{Proposition}[\bf solvability of subproblems in the graphical derivative-based Newton algorithm]\label{solv1} Let $\ph\colon\R^n\to\R$ be a ${\cal C}^{1,1}$ function around its tilt-stable local minimizer $\ox$. Then there exists a neighborhood $O$ of $\ox$ such that the set $\Upsilon(x)$ from \eqref{di1g} is nonempty and compact for all points $x\in O$.
\end{Proposition}\vspace*{-0.2in}
\begin{proof} It follows the lines in the proof of Proposition~\ref{exist} with the observation that
\begin{equation*}
(D\nabla\ph)^{-1}(v,x)(u)=Dg(v)(u)\;\mbox{ for all }\;u\in\R^n
\end{equation*}
in the notation of that proposition. Since the mapping $g$ therein is locally Lipschitzian, we derive from \cite[Proposition~9.24(a)]{rw} that its graphical derivative is nonempty-valued, closed-graph, and locally bounded. This yields the claimed properties of the sets $\Upsilon(x)$.
\end{proof}\vspace*{-0.07in}

One of the most important features of the graphical derivative-based Algorithm~\ref{nm2}, which distinguishes it from Algorithm~\ref{nm1}, is the possibility to supply the new algorithm with
an appropriate {\em subproblem} that resembles the one for the classical Newton method. To proceed, fix $x\in\R^n$ and consider the following optimization problem:
\begin{equation}\label{subp}
\min\;\;\ph(x)+\la v,w\ra+\sm\d^2\ph(x,v)(w)\;\; \mbox{subject to}\;\; w\in \R^n,
\end{equation}
where $v=\nabla\ph(x)$ and where $\d^2\ph(x,v)$ stands for the second subderivative of $\ph$ at $x$ for $v$ defined in Section~\ref{sect02}. If $w:=d$ is a stationary point of problem \eqref{subp}, then
$-v\in\sub \big(\frac{1}{2}\d^2\ph(x,v)\big)(d)$. Assuming further that $\ph$ is {\em twice epi-differentiable} at $x$ for $v$ and employing \cite[Theorem~13.40]{rw} tell us that $v\in(D\nabla\ph)(x)(d)$, and thus we get
\begin{equation*}
d\in(D\nabla\ph)^{-1}(v,x)(-v)=\Upsilon(x).
\end{equation*}
This indicates that a direction $d\in\R^n$ in Step~2 of Algorithm~\ref{nm2} can be calculated by solving the optimization problem \eqref{subp}. In the case where $\ph$ is a ${\cal C}^2$-smooth function, we have $\d^2\ph(x,v)(w)=\la\nabla^2\ph(x)w,w\ra$, which shows that subproblem \eqref{subp} reduces to the one in the classical Newton method for solving unconstrained optimization problems with ${\cal C}^2$-smooth objectives. Note that every directions from Algorithm~\ref{nm2} must be a stationary point of subproblem \eqref{subp} under the twice epi-differentiability of $\ph$, since in this case we always have by \cite[Theorem~13.40]{rw} that
\begin{equation*}
\big(D\nabla\ph\big)(x)(w)=\sub\Big(\frac{1}{2}\d^2\ph(x,v)\Big)(w),\;\;v=\nabla\ph(x).
\end{equation*}
Furthermore, under the assumptions of Theorem~\ref{solv2} below, the objective function of subproblem \eqref{subp} is strongly convex for all $x$ sufficiently close to $\ox$
as proven in the proof of Theorem ~\ref{solv2} and so the subproblem admits a unique optimal solution.
Since any direction of Algorithm~\ref{nm2} is a stationary point of subproblem \eqref{subp}, and since the objective function of it is convex, every direction is an optimal solution to \eqref{subp}, and thus it is unique.\vspace*{0.03in}

The next theorem tells us precisely that subproblem \eqref{subp} admits a {\em unique optimal solution} for all $x$ sufficiently close to the {\em tilt-stable} local minimizer of the function $\ph$ in question. Our proof below exploits a certain local monotonicity property closely related to tilt stability.

Recall that a mapping $T\colon\R^n\tto\R^m$ is {\em locally strongly monotone} with modulus $\tau>0$ around $(\ox,\oy)\in\gph T$ if there exist neighborhoods $U$ of $\ox$ and $V$ of $\oy$ such that
\begin{equation*}
\la x-u,v-w\ra\ge\tau\|x-u\|^2\;\mbox{ for all }\;(x,v),(u,w)\in(\gph T)\cap(U\times V).
\end{equation*}

Now we are ready to establish the aforementioned existence and uniqueness theorem for \eqref{subp} near $\ox$ and thus justify the possibility to find tilt-stable local minimizers of the original problem \eqref{main} by solving the much easier subproblem \eqref{subp} at each step of iterations.\vspace*{-0.07in}

\begin{Theorem}[\bf unique solutions of subproblems]\label{solv2} Let $\ph\colon\R^n\to\oR$ be a ${\cal C}^{1,1}$ function on a neighborhood of $\ox$, where $\ox$ is a tilt-stable local minimizer of $\ph$ with modulus $\kappa>0$, and let $\dd>0$ be such that for every $(x,v)\in(\gph\nabla\ph)\cap\B_\dd(\ox,0)$ the function $\ph$ is twice epi-differentiable at $x$ for $v=\nabla\ph(x)$.
Then for any $x$ close to $\ox$ subproblem \eqref{subp} admits a unique optimal solution.
\end{Theorem}\vspace*{-0.2in}
\begin{proof} Observe that the second-order growth characterization of tilt stability from Proposition~\ref{usogc2} implies by \eqref{ugc} that the gradient mapping $\nabla\ph$ is is locally strongly monotone with modulus $\kappa^{-1}$ around $(\ox,\ov)$ with $\ov=0$. Shrinking the neighborhoods $U$ and $V$ if necessary, suppose that $\ph$ is twice epi-differentiable at $x$ for $v$ whenever $(x,v)\in(\gph\nabla\ph)\cap(U\times V)$. Then employing \cite[Corollary~6.3]{pr96} tells us that for any such a pair $(x,v)$ the second subderivative $\frac{1}{2}\d^2\ph(x,v)$ is strongly convex
with modulus $\sm\kappa^{-1}$. Since functions $\ph$ of class ${\cal C}^{1,1}$ are continuously prox-regular at $\ox$ for $\ov=0$, we find $\ve>0$ and $\rho\ge 0$ such that
\begin{equation}\label{lower-quad}
\ph(u)\ge\ph(x)+\la v,u-x\ra-\frac{\rho}{2}\|x-u\|^2\;\mbox{ if }\;u\in\B_\ve(\ox)\;\mbox{ and }\;(x,v)\in(\gph\nabla\ph)\cap\B_{\ve}(\ox,0).
\end{equation}
This ensures the existence of  positive number $\nu$ with $\B_\nu(\ox,0)\subset U\times V$ for which
\begin{equation*}
\d^2\ph(x,v)(w)\ge-\rho\|w\|^2\;\mbox{ whenever }\;(x,v)\in(\gph\nabla\ph)\cap\B_\nu(\ox,0)\;\mbox{ and }\;w\in\R^n;
\end{equation*}
this is proved below in Proposition~\ref{prope} for any continuously prox-regular function $\ph$. Since $w\mapsto\d^2\ph(x,v)(w)$ is a positive homogeneous function of degree $2$, the above inequality yields $\d^2\ph(x,v)(0)=0$ for every pair $(x,v)\in(\gph\nabla\ph)\cap\B_\nu(\ox,0)$. This tells us that for all such pairs the second subderivative $\d^2\ph(x,v)$ is a proper function. Pick now any $x\in\R^n$ and define the new function $\psi_x\colon\R^n\to\oR$ by the second-order expansion
\begin{equation*}
\psi_x(w):=\ph(x)+\la v,w\ra+\sm\d^2\ph(x,v)(w)\;\mbox{ with }\;v=\nabla\ph(x),\;w\in\R^n.
\end{equation*}
The Lipschitz continuity of $\nabla\ph$ around $\ox$ gives us positive numbers $\dd_1$ and $\dd_2$ such that $\B_{\dd_1}(\ox)\times\B_{\dd_2}(\ov)\subset\B_\nu(\ox,0)$
and $\nabla\ph(x)\in\B_{\dd_2}(\ov)$ for all $x\in\B_{\dd_1}(\ox)$. Picking $x\in\B_{\dd_1}(\ox)$ and remembering that $v=\nabla\ph(x)$, we have $(x,v)\in(\gph\nabla\ph)\cap\B_{\nu}(\ox,0)$,
and thus the function $\psi_x$ is proper, l.s.c., and strongly convex. Since such functions admit unique optimal solutions, it finally verifies the claim of the theorem for any $x\in\B_{\dd_1}(\ox)$.
\end{proof}\vspace*{-0.07in}

Theorem~\ref{solv2} extends the well-known result on subproblems associated to the classical Newton method for \eqref{main} with ${\cal C}^2$-smooth objectives. As mentioned in Section~\ref{intro}, tilt stability for this setting amounts to the positive-definiteness of the Hessian matrix $\nabla^2\ph(\ox)$. Since subproblem \eqref{subp} reduces in the classical framework to
\begin{equation*}
\min\;\;\ph(x)+\la\nabla\ph(x),w\ra+\sm\la\nabla^2\ph(x)w,w\ra\;\; \mbox{subject to}\;\; w\in \R^n,
\end{equation*}
where $\nabla^2\ph(\ox)$ is positive-definite, the objective function of the latter subproblem is clearly strongly convex. Theorem~\ref{solv2} shows that a similar result is achieved for ${\cal C}^{1,1}$ functions if the second derivative is replaced by the second subderivative.\vspace*{0.03in}

Note that the solvability of subproblems in  Algorithm~\ref{nm2} was addressed in Proposition~\ref{exist} under the tilt stability assumption, while Theorem~\ref{solv2} goes much further in this vein. Indeed, it justifies a constructive way to find a required direction by solving subproblem \eqref{nm2} under an additional epi-differentiability assumption. We'll see in Section~\ref{sect05a} that the latter assumption holds automatically for a broad class of constrained optimization problems.\vspace*{0.03in}

Next we verify superlinear convergence of Algorithm~\ref{nm2} to tilt-stable local minimizers of ${\cal C}^{1,1}$ functions $\ph$ under the semismooth$^*$ property of the gradient mappings $\nabla\ph$. \vspace*{-0.07in}

\begin{Theorem}[\bf superlinear convergence of the graphical derivative-based Newton algorithm for ${\cal C}^{1,1}$ functions]\label{des10}
Let $\ph\colon\R^n\to\oR$ be a ${\cal C}^{1,1}$ function on a neighborhood of its tilt-stable local minimizer $\ox$ with modulus $\kappa>0$, and let $\nabla\ph$ be semismooth$^*$ at $\ox$. Then there exists $\dd>0$ such that for any starting point $x_0\in\B_\dd(\ox)$ we have that every sequence $\{x_k\}$ constructed by Algorithm~{\rm \ref{nm2}} converges to $\ox$ and the rate of  convergence is superlinear.
\end{Theorem}\vspace*{-0.2in}
\begin{proof} Since $\nabla\ph$ is semismooth$^*$ at $\ox$, by Proposition~\ref{semi5}(ii) we find $\dd'>0$ with
\begin{equation}\label{mss}
\|\nabla\ph(x)-\nabla\ph(\ox)-A_k(x-\ox)\|<\frac{1}{2\kappa}\|x-\ox\|\;\mbox{ for all }\;x\in\B_{\dd'}(\ox)\;\mbox{ and }\;A\in\co\Bar\nabla^2\ph(x).
\end{equation}
Pick $\dd>0$ from Theorem~\ref{inv02}(ii) for which estimate \eqref{gdi3} holds. Then Proposition~\ref{exist} gives us a neighborhood $O$ of $\ox$ with $\Upsilon(x)\ne\emp$ for all $x\in O$. Since $\nabla\ph$ is Lipschitz continuous around $\ox$, there are numbers $\dd_1\in(0,\dd')$ and $\dd_2>0$ such that $\B_{\dd_1}(\ox)\times\B_{\dd_2}(\ov)\subset\B_\dd(\ox,\ov)$, $\B_{\dd_1}(\ox)\subset O$,
and $\nabla\ph(x)\in\B_{\dd_2}(\ov)$ whenever $x\in\B_{\dd_1}(\ox)$. Letting $x_k\in\B_{\dd_1}(\ox)$ with $k\in\N$, we conclude from Proposition~\ref{solv1} that $\Upsilon(x_k)$ is nonempty. Pick $d_k\in\Upsilon(x_k)$ by Algorithm~{\rm\ref{nm1}} and set $x_{k+1}:=x_k+d_k$. Denoting $v_k:=\nabla\ph(x_k)$ tells us that $(x_k,v_k)\in(\gph\nabla\ph)\cap\B_\dd(\ox,\ov)$. By $\ov=0$ we get $d_k\in(D\nabla\ph)^{-1}(v_k,x_k)(\ov-v_k)$. It follows from estimate \eqref{gdi3} that for every   $k$ sufficiently large, there exists a matrix $A_k\in\co\Bar\nabla^2\ph(x_k)$ such that
\begin{equation}\label{newt}
\|x_k-\ox+d_k\|\le\kappa\,\|\nabla\ph(x_k)-\nabla\ph(\ox)-A_k(x_k-\ox)\|.
\end{equation}
Since $x_k\in\B_{\dd_1}(\ox)\subset\B_{\dd'}(\ox)$ for such $k$, we deduce from \eqref{mss} that
\begin{equation*}
\|x_{k+1}-\ox\|<\sm\|x_k-\ox\|,\;\mbox{ and so }\;x_{k+1}\in \B_{\dd_1}(\ox).
\end{equation*}
This tells us that every sequence $\{x_k\}$, generated by Algorithm~\ref{nm2} with the starting point $x_0\in\B_{\dd_1}(\ox)$, is contained in $\B_{\dd_1}(\ox)$ and converges to $\ox$ as $k\to\infty$. Combining it with Proposition~\ref{semi5}(ii) and estimate \eqref{newt} ensures that
\begin{equation*}
\|x_{k+1}-\ox\|=o(\|x_k-\ox\|)\;\mbox{ as }\;k\to\infty,
\end{equation*}
which verifies that the convergence is superlinear.
\end{proof}\vspace*{-0.17in}

\begin{Remark}[\bf comparison with related Newton-type algorithms]\label{newton-comp}{\rm Observe the following:

{\bf(i)} \cite[Theorem~13.57]{rw} tells us that if $\ph\colon\R^n\to\R$ is ${\cal C}^{1,1}$ around $x$ and twice epi-differentiable at $x$ for $\ov=\nabla\ph(x)$, then we have the inclusion
\begin{equation*}
(D\nabla\ph)(x)(w)\subset(D^*\nabla\ph)(x)(w)\;\mbox{ whenever }\;w\in\R^n.
\end{equation*}
This implies that $-\Upsilon(x)\subset\Upsilon^*(x)$ for such functions, and thus Algorithm~\ref{nm2} operates with a smaller set of directions in comparison with Algorithm~\ref{nm1} under the additional twice epi-differentiability assumption. Note to this end that the coderivative calculus in Algorithm~\ref{nm1} is much more developed than the one for the graphical derivatives. Observe also that, similarly to $\Upsilon^*$, the mapping $\Upsilon$ is neither a Newton approximation \cite{fp}, nor a Newton map \cite{kk}.

{\bf(ii)} The $B$-differential Newton method for solving equations developed in \cite{p90} is based on the $B$-derivative \cite{r}, which is actually the semiderivative in the terminology of \cite{rw} adapted here. As mentioned above, the latter construction reduces  for Lipschitzian mappings to the classical directional derivative. Thus, for unconstrained minimization problems with ${\cal C}^{1,1}$ objectives, the $B$-differential Newton method reduces to Algorithm~\ref{nm2}. The imposed assumptions in \cite{p90}, ensuring the existence of directions for each iteration, are rather restrictive and require that the mapping in question be actually strictly differentiable. The subsequent improvement in \cite{qi93} employs the so-called BD regularity assumption to achieve superlinear convergence. However, as  pointed out in \cite{qi93}, the imposed BD regularity assumption does not guarantee the solvability of subproblems in the proposed algorithm.

{\bf(iii)} The graphical/contingent derivative is also listed among generalized derivative constructions used in the Newton scheme developed in \cite{kk} to solve Lipschitzian equations, where some local convergence results are obtained under a set of assumptions not involving the fundamental notion of metric regularity in variational analysis; see, e.g., \cite{m06,rw}. The latter assumption is essentially used in \cite{hkmp} to ensure the solvability of subproblems in the generalized Newton method for equations that is based on the graphical derivative and provides superlinear convergence under assumptions different from \cite{kk} and Algorithm~\ref{nm2}. Extensions of \cite{hkmp} to solving set-valued inclusions are given in \cite{ds}. Note also a broad usage of metric regularity in Newton-type methods for Robinson's generalized equations; see \cite{dr}.}
\end{Remark}\vspace*{-0.05in}

Next we provide two examples of important classes of optimization problems that shed more light on how passing to subproblem \eqref{subp} brings us to a significantly simpler problem to solve.

\begin{Example}[\bf extended linear-quadratic programming] {\rm Define $\ph\colon\R^n\to\oR$ by
\begin{equation}\label{enlp}
\ph(x):=\la q,x\ra+\sm\la Qx,x\ra+f_{C,B}(b-Ax)\;\mbox{ for }\;x\in\R^n
\end{equation}
with the function $f_{{C,B}}\colon\R^m\to\oR$ given by
\begin{equation}\label{theta}
f_{{C,B}}(z):=\underset{p\in C}\sup{\big\{\la z,p\ra-\sm\la p,Bp\ra\big\}}\;\mbox{ for }\;z\in\R^m,
\end{equation}
where $q\in\R^n$, $b\in\R^m$, $Q$ is an $n\times n$ symmetric matrix, $A$ is an $n\times m$ matrix, $C$ is a polyhedral convex set in $\R^m$, and $B$ is an $m\times m$ symmetric and positive-semidefinite matrix. Problem \eqref{main} with the cost function \eqref{enlp} belongs to the class of {\em extended linear-quadratic programs} introduced by Rockafellar and Wets \cite{rw86}. Assuming in addition that $B$ is positive-definite, we deduce from the proof of \cite[Theorem~4.5]{mr} that $f_{C,B}$ is a ${\cal C}^{1,1}$ function.

Let $x\in\R^n$ and set $z:=b-Ax$, $u:=\nabla f_{C,B}(b-Ax)$, and $v:=\nabla\ph(x)$. Employing \cite[Theorem~5.4]{ms20} and \cite[Example~13.23]{rw}, we conclude for any $w\in\R^n$ that
\begin{eqnarray*}
\d^2\ph(x,v)(w)&=& \la Q w,w\ra+\d^2f_{C,B}(z,u)(-A w)\\
&=&\la Q w,w\ra+2f_{K_C(u,z-Bu),B}(-A w),
\end{eqnarray*}
where $K_C(u,z-B u)$ stands for the critical cone to $C$ at $u$ for $z-Bu$ defined by
\begin{equation*}
K_C(u,z-Bu)=T_C(u)\cap\{z-Bu\}^\perp.
\end{equation*}
This tells us that subproblems \eqref{subp} for \eqref{main} with the function $\ph$ from \eqref{enlp} can be simplified as
\begin{equation*}
\min\;\;\ph(x)+\la\nabla\ph(x),w\ra+\sm\la Qw,w\ra+f_{K_C(u,z-Bu),B}(-Aw)\;\;\mbox{subject to}\;\;w\in\R^n.
\end{equation*}
As seen, the main difference between the original problem \eqref{main} and subproblem \eqref{subp} is that the function $f_{{C,B}}$ in the objective $\ph$ is replaced by $f_{K_C(u,z-Bu),B}$. This means that the polyhedral set $C$ in the definition of $f_{{C,B}}$ is replaced with its critical cone, which often has a simpler structure than $C$.}
\end{Example}

\begin{Example}[\bf Augmented Lagrangians of constrained optimization problems] {\rm Recall that the {\em augmented Lagrangian} of the constrained optimization problem  \eqref{coop} with  $\psi\colon\R^n\to\R$ and $f\colon\R^n\to\R^m$ being  ${\cal C}^2$-smooth is defined by
\begin{equation}\label{aL}
\L(x,\lambda,\rho):=\psi(x)+e_{1/\rho}\dd_\Th(f(x)+\rho^{-1}\lambda)-\sm \rho^{-1}\Vert\lambda\Vert^2,
\end{equation}
where $(x,\lm,\rho)\in\R^n\times \R^m \times(0,\infty)$ and where $e_{1/\rho}\dd_\Th$ stands for the {\em Moreau envelope} of the indicator function $\dd_\Th$, which is defined
in \eqref{moreau} and broadly utilized in the following two sections. It is known that the augmented Lagrangian $\L$ is a ${\cal C}^{1,1}$ function with respect to $x$. Define the {\em Lagrangian} of \eqref{coop} by $L(x,\lm):=\psi(x)+\la\lm,f(x)\ra$ for any  $(x,\lm)\in\R^n\times\R^m$. To make our presentation easier, we assume further that $\Th$ is a polyhedral convex set but a similar observation holds for any parabolically regular set considered in Section~\ref{sect05a}. Remember also that for a given pair $(\lm,\rho)$ each iteration of the augmented Lagrangian method demands solving the problem
\begin{equation}\label{augm}
\min\;\;\ph(x):=\L(x,\lambda,\rho)\quad\mbox{subject to}\quad x\in\R^n.
\end{equation}
Set $v:=\nabla_x \L(x,\lambda,\rho)$ and $\mu:=\nabla \big(e_{1/\rho}\dd_\Th\big)\big(f(x)+\rho^{-1}\lambda\big)$. Employing \cite[Theorem~8.3]{mms19} shows that for any $w\in\R^n$ we get
\begin{eqnarray*}
\d^2_x\L\big((x,\lambda,\rho),v\big)&=&\big\langle\nabla^2_{xx}L(x,\mu)w,w\big\rangle+e_{1/{2\rho}}\big(\d^2\dd_\Th\big(f(x)+\rho^{-1}(\lm-\mu),\mu\big)\big)\big(\nabla f(x)w\big)\\
&=&\big\langle\nabla^2_{xx}L(x,\mu)w,w\big\rangle+ e_{1/{2\rho}}\big(\dd_{K_\Th (f(x)+\rho^{-1}(\lm-\mu),\mu)}\big)\big(\nabla f(x)w\big),
\end{eqnarray*}
where the last equality comes from \cite[Example~3.4]{mms19}, and where
\begin{equation*}
K_\Th\big(f(x)+\rho^{-1}(\lm-\mu),\mu\big)=T_\Th\big(f(x)+\rho^{-1}(\lm-\mu)\big)\cap\big\{\mu\big\}^\perp
\end{equation*}
is the critical cone to $\Th$ at $f(x)+\rho^{-1}(\lm-\mu)$ for $\mu$. Combining the above, we obtain the following representation for subproblem \eqref{subp} associated with $\ph$ that is taken from in \eqref{augm}:
\begin{equation*}
\underset{w\in\R^n}{\min}\;\;\ph(x)+\la\nabla\ph(x),w\ra+\sm\big\langle\nabla^2_{xx}L(x,\mu)w,w\big\rangle+e_{1/{\rho}}\big(\dd_{K_\Th(f(x)+\rho^{-1}(\lm-\mu),\mu)}\big)\big(\nabla f(x)w\big).
\end{equation*}
Thus comparing this subproblem with the original problem \eqref{augm} tells us that not only the terms $\psi(x)$ and $f(x)$ be replaced with the linear terms $\nabla^2\psi(x)w$ and $\nabla f(x)w$, but also the indicator function $\dd_\Th$ is substituted with $\dd_{K_\Th(f(x)+\rho^{-1}(\lm-\mu),\mu)}$. Note that the cone $K_\Th(f(x)+\rho^{-1}(\lm-\mu),\mu)$ often acquires a simpler structure in comparison with the original set $\Th$ in \eqref{aL}. Since a similar observation can be made for other important instances of $\Th$ such as the second-order cone defined by
\begin{equation*}
\Th={\cal Q}:=\big\{x=(y,x_n)\in\R^{n-1}\times\R\;\big|\;\|y\|\le x_n\big\},
\end{equation*}
we are going to demonstrate that what the transition from $\Th$ to its critical cone looks like for this set. Note that the second subderivative of $\dd_{\cal Q}$ is computed in  \cite[Example~5.8]{mms19}. Given  $(x,\lm)\in\R^n\times\R^m$, denote $\mu:=\nabla\big(e_{1/\rho}\dd_{\cal Q}\big)\big(f(x)+\rho^{-1}\lambda\big)$, which yields $\mu\in N_{\cal Q}\big(f(x)+\rho^{-1}(\lambda-\mu)\big)$. In order to see what simplifications can be provided for our subproblems in this case, we need to consider the following cases:

{\bf(a)} $f(x)+\rho^{-1}(\lambda-\mu)=0$. This tells us that $\mu\in-{\cal Q}$. If $\mu\in[\bd (-\Q )]\setminus\{0\}$, then
\begin{equation*}
K_\Q\big(f(x)+\rho^{-1}(\lm-\mu),\mu\big)=\big\{t(\mu',-\mu_m)\in \R^{m-1}\times \R|\, t\ge 0\big\}\;\;\mbox{with}\; \mu=(\mu',\mu_m),
\end{equation*}
which tells us that the cone ${\cal Q}$ is replaced with a ray in our subproblem. If $\mu\in\inte(-\Q)$, we obtain
\begin{equation*}
K_\Q\big(f(x)+\rho^{-1}(\lm-\mu),\mu\big)=\{\mu\}^\perp,
\end{equation*}
and thus $\Q$ is replaced with a hyperplane in our subproblem. If $\mu=0$, then$K_\Th\big(f(x)+\rho^{-1}(\lm-\mu),\mu\big)=\Q$, and so no change occurs.

{\bf(b)} $f(x)+\rho^{-1}(\lambda-\mu)\in\big(\bd\Q\big)\setminus\{0\}$. It is not hard to see in this case that
\begin{equation*}
K_\Q\big(f(x)+\rho^{-1}(\lm-\mu),\mu\big)=
\begin{cases}
\{\mu\}^\perp&\mbox{if}\;\;\mu\ne 0,\\
T_\Q\big(f(x)+\rho^{-1}(\lm-\mu)\big)&\mbox{if}\;\;\mu=0.
\end{cases}
\end{equation*}
This shows that $\Q$ is replaced with either a hyperplane (when $\mu\ne 0$) or a closed halfspace (when $\mu=0$) in our subproblem.

{\bf(c)} $f(x)+\rho^{-1}(\lambda-\mu)\in\inte \Q$. This implies that $\mu=0$, and so we arrive at
\begin{equation*}
K_\Q\big(f(x)+\rho^{-1}(\lm-\mu),\mu\big)=\R^m.
\end{equation*}}
\end{Example}\vspace*{0.01in}

To conclude this section, we should mention that some {\em globalization} strategies can be combined with the local superlinear convergence established in Theorem~\ref{des10} to derive global   convergence of  Algorithm~\ref{nm2}. A common method here is to use a {\em line search} strategy and update the sequence $\{x_k\}$ from Theorem~\ref{des10} by $x_{k+1}:=x_k+\al_kd_k$, where $\al_k$ is a stepsize in the direction $d_k$. One of the most popular line search is to choose a stepsize $\al_k$ that satisfies the condition
\begin{equation}\label{armijo}
\ph(x_{k}+ \al_kd_k)\le\ph(x_k)+\mu\al_k\la\nabla\ph(x_k),d_k\ra,
\end{equation}
where $\mu\in(0,1)$. This condition is referred to as the {\em Armijo rule}. To derive global superlinear convergence, it is required to show that in a neighborhood of the limit point of the sequence $\{x_k\}$ from Algorithm~{\rm\ref{nm2}}, which is $\ox$ under the assumptions utilized in Theorem~\ref{des10}, the unit stepsize will be accepted in \eqref{armijo}; namely, this estimate holds for $\al_k=1$ for all $k$ sufficiently large. The final result of this section aims at furnishing such a conclusion for the latter sequence.
\vspace*{-0.05in}

\begin{Proposition}[\bf acceptance of unit stepsize in the Armijo rule]\label{desc} Let $\ph\colon\R^n\to\oR$ be a ${\cal C}^{1,1}$ function with Lipschitz constant $\ell>0$ for $\nabla\ph$ around its tilt-stable local minimizer $\ox$ with modulus $\kappa>0$, and let $\nabla\ph$ be semismooth$^*$ at $\ox$. Assume further that the sequence $\{x_k\}$ is generated by Algorithm~{\rm\ref{nm2}} with $x_k\ne \ox$ for all $k\in\N$. Then for any $\mu\in\big(0,1/(4\ell\kappa)\big)$ there exists $\bar k\in\N$ such that whenever $k\ge\bar k$ we have
\begin{equation*}
\ph(x_{k}+ d_k)\le\ph(x_k)+\mu\la\nabla\ph(x_k),d_k\ra,
\end{equation*}
where the directions $d_k$ are taken from Algorithm~{\rm\ref{nm2}}.
\end{Proposition}\vspace*{-0.17in}
\begin{proof} Using Theorem~\ref{des10}, we find $\dd>0$ such that for any $x_0\in\B_\dd(\ox)$ every sequence $\{x_k\}$ constructed by Algorithm~\ref{nm2} satisfies the condition
\begin{equation*}
\|x_{k+1}-\ox\|=o(\|x_k-\ox\|)\;\;\mbox{as}\;\;k\to\infty.
\end{equation*}
This implies that $\lim_{k\to\infty}{\|d_k\|}/{\|x_k-\ox\|}=1$, and hence $\lim_{k\to\infty}{\|x_k-\ox\|}/{\|d_k\|}=1$. This  yields
\begin{equation}\label{qam}
\|x_{k+1}-\ox\|=o(\|d_k\|)\;\;\mbox{as}\;\;k\to \infty
\end{equation}
and $|\frac{\|x_k-\ox\|}{\|d_k\|}-1|\le{{1}/{(4\ell\kappa)}}$ for all $k$ sufficiently large. Let $\mu\in(0,{1}/{(4\ell\kappa)})$ and observe that $\ell\kappa\ge 1$. Thus for all $k$ sufficiently large we arrive at
\begin{equation}\label{ed02}
\mu\frac{\|x_k-\ox\|}{\|d_k\|}-{\frac{1}{2\ell\kappa}}<\frac{\mu}{4\ell\kappa}+\mu-\frac{1}{2\ell\kappa}\le\frac{1}{4\ell\kappa}+\mu-\frac{1}{2\ell\kappa}=\mu-\frac{1}{4\ell\kappa}<0.
\end{equation}
Assume without loss of generality that  $\B_{\dd}(\ox)\times\B_{\dd\ell}(\ov)\subset U\times V$, where $\ov:=0$ and the neighborhoods $U$ and $V$ are taken from Proposition~\ref{usogc2}.
Since the sequence $\{x_k\}$ converges to  $\ox$, we obtain $(x_k,\nabla \ph(x_k))\in \B_{\dd}(\ox)\times\B_{\dd\ell}(\ov)$ for all $k$ sufficiently large.  Since $x_{k+1}=x_k+d_k$,  the uniform second-order growth condition from Proposition~\ref{usogc2} ensures that
\begin{eqnarray*}
\ph(x_{k+1})-\ph(x_k)-\mu\la \nabla\ph(x_k),d_k\ra&\le&\la \nabla\ph(x_{k+1}),d_k\ra-\frac{1}{2\kappa}\|d_k\|^2- \mu\la\nabla\ph(x_k),d_k\ra\\
&\le&\|\nabla\ph(x_{k+1})\|\cdot\|d_k\|-\frac{1}{2\kappa}\|d_k\|^2+\mu\|\nabla\ph(x_k)\|\cdot\|d_k\|\\
&\le&\ell\|x_{k+1}-\ox\|\cdot\|d_k\|-\frac{1}{2\kappa}\|d_k\|^2+\mu\ell\|x_k-\ox\|\cdot\|d_k\|\\
&=& \frac{o(\|d_k\|^2)}{\|d_k\|^2}+\ell\Big(\mu\frac{\|x_k-\ox\|}{\|d_k\|}-{\frac{1}{2\ell\kappa}}\Big),
\end{eqnarray*}
where the last equality results from \eqref{qam}. Combining this with \eqref{ed02} justifies the claimed estimate for all $k$ sufficiently large and thus completes the proof.
\end{proof}\vspace*{-0.2in}

\section{Newton Algorithms for Prox-Regular Functions}\sce\label{sect05}\vspace*{-0.05in}

In this section we proceed with extensions of both Algorithms~\ref{nm1} and \ref{nm2} to a much more general class of continuously prox-regular functions $\ph\colon\R^n\to\oR$. This framework encompasses problems of constrained optimization, which are explicitly considered in the next section. Such an extension is based on the remarkable facts of variational analysis allowing us to pass from objective functions of class ${\cal C}^{1,1}$ to continuously prox-regular objectives by using Moreau envelopes. Recall that the {\em Moreau envelope} of $\ph\colon\R^n\to\oR$ for $r>0$ is defined by the infimal convolution
\begin{equation}\label{moreau}
e_r\ph(x):=\inf_w\Big\{\ph(w)+\frac{1}{2r}\|w-x\|^2\Big\},
\end{equation}
and the corresponding {\em proximal mapping} of $\ph$ is given by
\begin{equation}\label{proxmap}
P_r\ph(x):=\mbox{argmin}_w\Big\{\ph(w)+\frac{1}{2r}\|w-x\|^2\Big\}.
\end{equation}

The following result, which is taken from \cite[Proposition~13.37]{rw}, collects the needed properties of Moreau envelopes used below.\vspace*{-0.07in}

\begin{Proposition}[\bf Moreau envelopes for prox-regular functions]\label{moru} Let $\ph\colon\R^n\to\oR$ be continuously prox-regular at $\ox$ for $\ov=0$, and let $\ph$ be bounded from below by a quadratic function on $\R^n$. Then for any $r>0$ sufficiently small there exists an $r$-dependent neighborhood $U$ of $\ox$ on which $e_r\ph$ is of class ${\cal C}^{1,1}$ with $\nabla e_r\ph(\ox+r\ov)=\ov$, and we have the representation
\begin{equation}\label{eq02}
\nabla e_r\ph(u)=\big(r I+T^{-1}\big)^{-1}(u)\;\mbox{ for all }\;u\in U,
\end{equation}
where $T$ is a graphical localization of $\sub\ph$ around $(\ox,\ov)$.
\end{Proposition}\vspace*{-0.07in}

Observe that the boundedness from below assumption on $\ph$ in Proposition~\ref{moru} is not restrictive and will be dropped in this section. Indeed, since our analysis depends only on the local geometry of $\gph\sub\ph$ around $(\ox,\ov)$, by adding to $\ph$ the indicator of some compact neighborhood of $\ox$ if necessary, we can make $\ph$ to be bounded from below by a quadratic function on $\R^n$. In what follows we always assume that there is $\rho\ge 0$ with
\begin{equation*}
\ph(x)\ge\ph(\ox)-\frac{\rho}{2}\|x-\ox\|^2\;\mbox{ for all }\;x\in\R^n.
\end{equation*}

Thus the usage of the Moreau envelope \eqref{moreau} allows us to pass from the original problem \eqref{main} with a continuously prox-regular objective to the similarly formulated problem:
\begin{equation}\label{main2}
\mbox{minimize }\;e_r\ph(x)\;\mbox{ subject to }\;x\in\R^n
\end{equation}
with the objective given by a ${\cal C}^{1,1}$ function. Let us emphasize again that, although both problems \eqref{main} and \eqref{main2} are written in the same unconstrained optimization format, they are significantly different from each other due to the actual constrained and highly nonsmooth nature of \eqref{main} in the case of continuously prox-regular objectives. To proceed with the implementation of Algorithms~\ref{nm1} and \ref{nm2} for problem \eqref{main} via the passage to \eqref{main2}, we have to find appropriate assumptions on $\ph$ in \eqref{main}, which ensure the fulfillment of those in \eqref{main2} allowing us to apply the results of Sections~\ref{sect04} and \ref{sect04a} to \eqref{main2}. It luckily occurs that the corresponding assumptions are the same. This is shown in the proof of the following major result. Note that, similarly to Section~\ref{sect04a}, we can supply Algorithm~\ref{nm2} for prox-regular functions regularized via Moreau envelopes by the corresponding {\em subproblem} of type \eqref{subp} under an additional {\em twice epi-differentiability} assumption. We'll proceed in more details in this direction in Section~\ref{sect05a} for constrained optimization.\vspace*{-0.07in}

\begin{Theorem}[\bf solvability and superlinear convergence of Newton algorithms for prox-regular functions]\label{supunc} Let $\ph\colon\R^n\to\oR$ be continuously prox-regular at $\ox$ for $\ov=0$, where $\ox$ is a tilt-stable local minimizer for $\ph$ with modulus $\kappa>0$. Assume further that the mapping $\sub\ph$ is semismooth$^*$ at $(\ox,\ov)$. Then given any $r>0$ sufficiently small, there exists $\dd>0$ such that for each starting point $x_0\in\B_\dd(\ox)$ both Algorithms~{\rm\ref{nm1}} and {\rm\ref{nm2}} for \eqref{main2} are well-defined, and every sequence $\{x_k\}$ constructed by either of them for the function $e_r\ph$ converges to $\ox$ and the rate of convergence is superlinear.
\end{Theorem}\vspace*{-0.17in}
\begin{proof}
We split the proof of the theorem into several steps, which are formulated as claims of their own interest, and begin by showing that the property of tilt stability is disseminated from a continuous prox-regular function to its Moreau envelope.\\[0.03in]
{\bf Claim~1:} {\em If $\ph$ is continuously prox-regular at $\ox$ for $\ov=0$ and if $\ox$ is a tilt-stable local minimizer for $\ph$ with modulus $\kappa>0$, then for any $r>0$ sufficiently small the point $\ox$ is a tilt-stable local minimizer for $e_r\ph$ with modulus $\kappa+2r$}.\\[0.03in]
To verify this claim, pick any small $r>0$ from Proposition~\ref{moru} and by representation \eqref{eq02} find neighborhoods $U$ of $\ox+r\ov=\ox$ and $V$ of $\ov=0$ such that for all $(x,v)\in U\times V$ we have
\begin{equation}\label{eq03}
v=\nabla(e_r\ph)(x)\iff v\in\sub\ph(x-rv).
\end{equation}
It follows from the Fermat rule $0\in\sub\ph(\ox)$ that $\nabla(e_r\ph)(\ox)=0$. Taking this into account, we deduce from the ${\cal C}^{1,1}$ property of \eqref{moreau} that $e_r\ph$ is continuously prox-regular at $\ox$ for $0$. Define further the linear transformation ${\cal L}\colon\R^n\times\R^n\to\R^n\times\R^n$ by ${\cal L}(x,v):=(x-r v,v)$. This allows us to equivalently rewrite \eqref{eq03} as
\begin{equation*}
\big(\gph\nabla(e_r\ph)\big)\cap\big(U\times V\big)=\big\{(x,v)\in U\times V\big|\;{\cal L}(x,v)\in\gph\sub\ph\big\}.
\end{equation*}
Since for any $(x,v)\in\R^n\times\R^n$ we obviously have
\begin{equation*}
\nabla{\cal L}(x,v)=\left(\begin{array}{cc}
I&-r I\\
0&I\\
\end{array}\right)
\end{equation*}
with $I$ standing for the $n\times n$ identity matrix, the Jacobian matrix $\nabla{\cal L}(x,v)$ is of full rank. Appealing to \cite[Exercise~6.7]{rw} tells us that
\begin{equation}\label{rd01}
T_{\scriptsize{\gph\nabla(e_r\ph)}}(x,v)=\big\{(w,u)\in\R^n\times\R^n\big|\;(w-ru,u)\in T_{\scriptsize{\gph\sub\ph}}(x-rv,v)\big\}
\end{equation}
for any $(x,v)\in U\times V$. Since $\ox$ is a tilt-stable local minimizer of $\ph$ with modulus $\kappa>0$, it follows from Proposition~\ref{tiltch}(ii) that
there exists $\dd>0$ such that
\begin{equation}\label{cso3}
\la w,q\ra\ge\frac{1}{\kappa}\|w\|^2\;\;\mbox{whenever}\;\;(w,q)\in T_{\scriptsize{\gph\sub\ph}}(u,v)\;\;\mbox{with}\;(u,v)\in\big(\gph\sub\ph\big)\cap\B_\dd(\ox,0).
\end{equation}
Shrinking the neighborhoods $U$ and $V$ if necessary, suppose without loss of generality that ${\cal L}(U\times V)\subset\B_\dd(\ox,\ov)$. Picking $(x,v)\in(\gph\nabla(e_r\ph))\cap(U\times V)$
and $(w,u)\in T_{\scriptsize{\gph\nabla(e_r\ph)}}(x,v)$, we deduce from \eqref{rd01} the inclusions
\begin{equation*}
(w-ru,u)\in T_{\scriptsize{\gph\sub\ph}}(x-rv,v)\;\mbox{ and }\;(x-r v,v)\in{\cal L}(U\times V)\subset\B_\dd(\ox,0).
\end{equation*}
Employing now (\ref{cso3}) brings us to the estimate
\begin{equation*}
\la u,w-ru\ra\ge\frac{1}{\kappa}\|w-ru\|^2,
\end{equation*}
which in turn implies that
\begin{equation*}
\la w,u\ra\ge\frac{1}{\kappa+2r}\|w\|^2.
\end{equation*}
In summary, we arrive at the inequality
\begin{equation*}
\la w,u\ra\ge\frac{1}{\kappa+2r}\|w\|^2\;\;\mbox{for all}\;\;(w,u)\in T_{\scriptsize{\gph\nabla(e_r\ph)}}(x,v)\;\;\mbox{with}\;(x,v)\in\big(\gph(\nabla e_r\ph)\big)\cap(U\times V),
\end{equation*}
which ensures by Proposition~\ref{tiltch} that $\ox$ is a tilt-stable local minimizer of the Moreau envelope $e_r\ph$ with modulus $\kappa+2r$. This verifies the claim.\vspace*{0.05in}

To proceed further with the proof of theorem, let us show that the required semismooth$^*$ property for \eqref{main} is equivalent to the same property for \eqref{main2}.\\[0.03in]
{\bf Claim~2:} {\em In the setting of the theorem we have that for any $r>0$ sufficiently small the semismooth$^*$ property of $\nabla(e_r\ph)$ at $\ox$ is equivalent to the semismooth$^*$ property of $\sub\ph$ at $(\ox,\ov)$.}\\[0.03in]
The proof of the claimed equivalence fully relies on representation \eqref{eq02}. Pick a small $r>0$ for which \eqref{eq02} is satisfied. If $\nabla(e_r\ph)$ is semismooth$^*$ at $\ox$, then the latter equality tells us that $rI+T^{-1}$ is semismooth$^*$ at $(\ov,\ox+r\ov)$. It follows from \cite[Proposition~3.6]{go19} that the mapping $T^{-1}$ is semismooth$^*$ at $(\ov,\ox)$ and so is $T$ at $(\ox,\ov)$. Since $T$ is a graphical localization of $\sub\ph$ around $(\ox,\ov)$, the semismoothness$^*$ of $T$ at $(\ox,\ov)$ is equivalent to that for $\sub\ph$ at the same point. This verifies the semismoothness$^*$ of $\sub\ph$ at $(\ox,\ov)$. The converse implication is justifies similarly.\vspace*{0.05in}

Now we are ready to finalize the proof of theorem.\\[0.03in]
{\bf Claim~3:} {\em Both Newton-type algorithms for \eqref{main2} are well-defined and superlinearly convergent.} Since we know from Claim~1 that $\ox$ is a tilt-stable local minimizer for $e_r\ph$ whenever $r>0$ is sufficiently small, the solvability of subproblems in Algorithm~\ref{nm1} and Algorithm~\ref{nm2} for problem \eqref{main2} follows from Proposition~\ref{exist} and Proposition~\ref{solv1}, respectively. Furthermore, Claim~2 tells us that the mapping $\nabla(e_r\ph)$ is semismooth$^*$ at $\ox$ under the assumptions of the theorem. Thus we deduce the asserted convergence and the rate of convergence of these algorithms for problem \eqref{main2} by applying the corresponding statements of Theorem~\ref{des12} or Theorem~\ref{des10}.
\end{proof}\vspace*{-0.05in}

As mentioned in Section~\ref{intro}, the regularization procedure of type \eqref{main2} was first suggested in \cite{fq} for finite-valued convex functions on the base of the semismooth Newton method. A drawback of that paper, in contrast to our developments above, is that the imposed assumptions in \eqref{main2} were not expressed in terms of the original problem \eqref{main}, but via the data of the regularized one \eqref{main2}. Observe also a constructive approach of \cite{fq} to find an approximate solution to the optimization problem in definition \eqref{moreau} of the Moreau envelope in the case of convex functions $\ph$. In our future research we intend to develop a similar approach to numerical implementations of Algorithm~\ref{nm1} and Algorithm~\ref{nm2} for continuously prox-regular functions. Finally, we refer the reader to \cite{hl} and the bibliography therein for other developments on the computation of Moreau envelopes for piecewise linear-quadratic functions and their conjugates.\vspace*{-0.15in}

\section{Applications to Constrained Optimization}\sce\label{sect05a}\vspace*{-0.05in}

Here we present some applications of the Newton-type method based on the graphical derivative, which was developed in Sections~\ref{sect04a} and \ref{sect05}, to solving {\em constrained optimization problems} defined in \eqref{coop}, where $\psi\colon\R^n\to\R$ and $f\colon\R^n\to\R^m$ are ${\cal C}^2$-smooth around the references points, and where $\Th\subset\R^m$ is closed and convex. As mentioned in Section~\ref{intro}, the case where $\Th$ is a cone refers to the class of {\em conic programs} highly important in optimization theory and applications; see, e.g., \cite{bs}. Problem \eqref{coop} can be obviously rewritten in the unconstrained form
\begin{equation}\label{coop2}
\mbox{minimize }\;\ph(x):=\psi(x)+\dd_\O(x)\;\mbox{ with }\;\O:=\big\{x\in\R^n\big|\;f(x)\in\Th\big\},
\end{equation}
where $\ph$ is continuously prox-regular at the points in question as follows from \cite[Section~13F]{rw}.

As shown in Section~\ref{sect04a}, a constructive realization of Algorithm~\ref{nm2} for ${\cal C}^{1,1}$ functions (and hence of its extension for prox-regular ones in Section~\ref{sect05}) requires an additional assumption on the {\em twice epi-differentiability} of the cost function. To obtain efficient conditions for the twice epi-differentiability of the function $\ph$ from \eqref{coop2} in terms of the initial data of this problem, we invoke our recent developments \cite{mms19} on {\em parabolic regularity} in second-order variational analysis.

Recall that a set $\Th\subset\R^m$ is {\em parabolically regular} at $\oz\in\Th$ for $\ou\in\R^m$ if whenever $w\in\R^m$ is such that $\d^2\dd_\Th(\bar z,\ou)(w)<\infty$ there exist, among the sequences $t_k\dn 0$ and $w_k\to w$ with $\Delta_{t_k}^2\dd_\Th(\bar z,\ou)(w_k)\to\d^2\dd_\Th(\bar z,\ou)(w)$, those with the additional property that $\limsup_{k\to\infty}\|w_k-w\|/t_k<\infty$ as $k\to\infty$.
This notion was introduced in \cite[Definition~13.65]{rw}, but was not explored in \cite{rw} or anywhere else (before the recent paper \cite{mms19}) further than in the fully amenable setting. The paper \cite{mms19} offers an extensive study of parabolic regularity by revealing its remarkable properties as an appropriate second-order regularity notion for a large class of sets that overwhelmingly appear in variational analysis and constrained optimization. This class strictly includes all the  ${\cal C}^2$-{\em cone reducible} sets in the sense of \cite[Definition~3.135]{bs} and encompasses convex polyhedra, the second-order cone, the cone of symmetric and positive semidefinite matrices, etc. Furthermore, parabolic regularity, combined with related developments of \cite{mms19a}, occurs to be very instrumental in the study and calculations of second subderivatives and twice epi-differentiability of functions while being employed in our numerical applications given below.\vspace*{0.03in}

Before formulating the needed results in this direction, let us recall some additional notions and constructions. Given $\O$ from \eqref{coop2}, pick $(x,v)\in\gph N_\O$ and define the set of {\em Lagrange multipliers} associated with $(x,v)$ by
\begin{equation*}
\Lambda(x,v):=\big\{\lambda\in N_\Th\big(f(x)\big)|\;\nabla f(x)^*\lambda=v\big\}.
\end{equation*}
Recall that a set $\Th\subset\R^m$ is {\em parabolically derivable} at $\oz$ for $w$ if $T_\Th^2(\oz,w)\ne\emp$ is for each $u\in T_\Th^2(\oz,w)$ there are $\ve>0$ and $\xi\colon[0,\ve]\to\Th$ with $\xi(0)=\oz$, $\xi'_+(0)=w$, and $\xi''_+(0)=u$, where
\begin{equation*}
\xi''_+(0):=\lim_{t\dn 0}\frac{\xi(t)-\xi(0)-t\xi'_+(0)}{\sm t^2},
\end{equation*}
and where $T_\Th^2(\oz,w)$ is the {\em second-order tangent set} to $\Th$ at $\oz$ for $w\in T_\Th(\oz)$ given by
\begin{equation*}
T_\Th^2(\oz,w):=\big\{u\in\R^m|\;\exists\,t_k{\downarrow}0,\;\;u_k\to u\;\;\mbox{ as }\;k\to\infty\;\;\mbox{with}\;\;\ox+t_kw+\sm t_k^2u_k\in\Th\big\}.
\end{equation*}
Parabolic derivability is a fairly common property in second-order analysis; see, e.g., \cite{rw}.

Recall also that a set-valued mapping $F\colon\R^n\tto\R^m$ is {\em metrically subregular} at $(\ox,\oy)\in\gph F$ if there exists a neighborhood $U$ of $\ox$ and a number $\mu>0$ such that
\begin{equation*}
{\rm dist}\big(x;F^{-1}(\oy)\big)\le\mu\,{\rm dist}\big(\oy;F(x)\big)\;\mbox{ for all }\;x\in U.
\end{equation*}

The following proposition collects some results from \cite[Theorem~5.6 and Corollary~5.11]{mms19} ensuring the twice epi-differentiability of parabolically regular constrained systems of the type in \eqref{coop2} and calculating their second subderivatives.\vspace*{-0.07in}

\begin{Proposition}[\bf twice epi-differentiability of constraint systems]\label{tedi} Let $\O$ be taken from \eqref{coop2}, and let $(\ox,\ov)\in\gph N_\O$.
Assume further that:

{\bf(a)} The set-valued mapping $x\mapsto f(x)-\Th$ is metrically subregular at $(\ox,0)$.

{\bf(b)} There exists a positive number $\ve$ such that for any $(x,v)\in(\gph N_\O)\cap\B_\ve(\ox,\ov)$ and any $\lm\in\Lambda(x,v)$ the set $\Th$ is parabolically regular at $f(x)$ for $\lm$ while being also parabolically derivable at $f(x)$ for each $w\in T_\Th(f(x))\cap\{\lm\}^\bot$.\\[0.5ex]
Then for any $(x,v)\in(\gph N_\O)\cap\B_\ve(\ox,\ov)$ the indicator function $\dd_\O$ is properly twice epi-differentiable at $x$ for $v$ and its second subderivative at $x$ for $v$ is calculated by
\begin{equation*}
\d^2\dd_\O(x,v)(w)=\max_{\lm\in\Lambda(x,v)}\big\{\langle\lm,\nabla^2f(x)(w,w)\rangle+\d^2\dd_\Th\big(f( x),\lm\big)\big(\nabla f(x)w\big)\big\},\quad w\in\R^n.
\end{equation*}
Furthermore, we have the domain representation $\dom\d^2\dd_\O(x,v)=T_\O(x)\cap\{v\}^\bot$.
\end{Proposition}\vspace*{-0.07in}

Remembering that $\psi$ in \eqref{coop2} is ${\cal C}^2$-smooth around $\ox$ and denoting the {\em Lagrangian} of \eqref{coop2} by $L(x,\lm):=\psi(x)+\la\lm,f(x)\ra$ as $(x,\lm)\in\R^n\times\R^m$, we deduce from Proposition~\ref{tedi} with $\ov:=-\nabla\psi(\ox)$ by employing the elementary sum rule for second subderivatives from \cite[Exercise~13.18]{rw} that there exists $\ve>0$ such that for every $(x,v)\in(\gph\sub\ph)\cap\B_\ve(\ox,0)$ the function $\ph$ is twice epi-differentiable at $x$ for $v$ and its second subderivative is calculated by
\begin{eqnarray}
\d^2\ph(x,v)(w)&=&\la\nabla^2\psi(x)w,w\ra+\max_{\lm\in\Lambda(x,v-\nabla\ph(x))}\big\{\la\lm,\nabla^2f(x)(w,w)\ra+\d^2\dd_\Th\big(f( x),\lm\big)\big(\nabla f(x)w\big)\big\}\nonumber\\
&=&\max_{\lm\in\Lambda(x,v-\nabla\psi(x))}\big\{\la\nabla_{xx}^2L(x,\lm)w,w\ra+\d^2\dd_\Th\big(f(x),\lm\big)\big(\nabla f(x)w\big)\big\},\quad w\in\R^n\label{ssf}.
\end{eqnarray}

Using the above discussions and applying Algorithm~\ref{nm2} to problem \eqref{coop2} regularized via \eqref{moreau}, with taking into account the direction search by Theorem~\ref{solv2}, we arrive at the following Newton-type algorithm for constrained optimization based on {\em second subderivatives}.\vspace*{-0.07in}

\begin{Algorithm}[\bf second subderivative-based Newton method for constrained optimization problems]\label{nm3} {\rm Considering problem \eqref{coop2} under the assumptions above, let $x_0\in\R^n$, set $k:=0$, and pick any $r>0$.\\
{\bf Step~1:} If $0\in\sub\ph(x_k)$, then stop.\\
{\bf Step~2:} Otherwise, let $v_k=\nabla(e_r\ph)(x_k)$, select $w_k$ as a stationary point of the subproblem
\begin{equation}\label{subp3}
\underset{w\in\R^n}{\min}\;\la v_k,w\ra+\sm\d^2\ph(x_k-rv_k,v_k)(w),
\end{equation}
and then set $d_k:=w_k-rv_k$ and $x_{k+1}:=x_k+d_k$.\\
{\bf Step~3:} Let $k\leftarrow k+1$ and then go to Step~1.}
\end{Algorithm}\vspace*{-0.07in}

Observe by \eqref{eq03} that the stationary condition $0\in\sub\ph(x_k)$ amounts to $\nabla(e_r\ph)(x_k)=0$, and that it is equivalently expressed via the initial data of \eqref{coop} by
\begin{equation*}
-\nabla\psi(x_k)\in\nabla f(x_k)^*N_\Th\big(f(x_k)\big)
\end{equation*}
under the qualification condition $N_\Th(f(x_k))\cap\ker\nabla f(x_k)^*=0$, which holds, in particular, when $\nabla f(x_k)$ is of full rank; see \cite[Corollary~3.13]{m18}. As we see from \eqref{ssf}, the second subderivative in \eqref{subp3} is calculated in terms of the initial data of \eqref{coop}. It also follows from \eqref{subp3} that if $w_k$ is a stationary point of this subproblem, then we get the inclusion
\begin{equation*}
-v_k\in\sub\big(\sm\d^2\ph(x_k-rv_k,v_k)\big)(w_k).
\end{equation*}
Furthermore, the assumptions of Proposition~\ref{tedi} ensures that the above inclusion can be equivalently rewritten by \cite[Theorem~13.40]{rw} as
\begin{equation*}
-v_k\in(D\sub\ph)(x_k-rv_k,v_k)(w_k)=(D\sub\ph)(x_k-rv_k,v_k)(d_k+rv_k),
\end{equation*}
where the direction $d_k$ is taken from Algorithm~\ref{nm3}. Employing now \eqref{rd01} tells us that the latter inclusion amounts to $-v_k\in(D\nabla(e_r\ph))(x_k)(d_k)$,
which confirms that Algorithm~\ref{nm3} is actually Algorithm~\ref{nm2} implemented for $e_r\ph$ with $\ph$ taken from \eqref{coop2}.\vspace*{0.03in}

Remember that in the proof of Theorem~\ref{solv2} we postponed the verification of the quadratic lower estimate \eqref{lower-quad} for the second subderivative of ${\cal C}^{1,1}$ functions. As promised, now we establish such an estimate for the general case of continuously prox-regular functions. The following proposition and its proof extend those given in \cite[Theorem~4.1]{mms19a}. The obtained result shows that subproblem \eqref{subp3} always admit an optimal solution under the tilt stability assumption.\vspace*{-0.05in}

\begin{Proposition}[\bf properties of second subderivatives of prox-regular functions]\label{prope} Let $\ph\colon\R^n\to\oR$ be continuously prox-regular at $\ox$ for $\ov\in\partial\ph(\ox)$. Then there exist $\ve>0$ and $\rho\ge 0$ such that for every $(x,v)\in(\gph\sub\ph)\cap\B_\ve(\ox,\ov)$ the second subderivative $\d^2\ph(x,v)$ is a proper and l.s.c. function satisfying the quadratic lower estimate
\begin{equation}\label{quad-est}
\d^2\ph(x,v)(w)\ge-\rho\|w\|^2\;\mbox{ whenever }\;w\in\R^n.
\end{equation}
\end{Proposition}\vspace*{-0.17in}
\begin{proof} The claimed lower semicontinuity of $\d^2\ph(x,v)$ follows from \cite[Proposition~13.5]{rw}. To verify the lower estimate \eqref{quad-est}, deduce from the assumed continuous prox-regularity of $\ph$ the existence of $\ve>0$ and $\rho\ge 0$ ensuring that
\begin{equation*}
\ph(u)\ge\ph(x)+\la v,u-x\ra-\frac{\rho}{2}\|x-u\|^2\;\mbox{ for all }\;u\in\B_\ve(\ox),\;(x,v)\in(\gph\sub\ph)\cap\B_{\ve}(\ox,\ov).
\end{equation*}
Picking $(x,v)\in(\gph\sub\ph)\cap\B_{\ve/2}(\ox,\ov)$ and $w\in \R^n$, deduce from the above inequality that whenever $t>0$ is sufficiently small and $w'$ is close to $w$ we get
\begin{equation*}
\Delta_t^2\ph(x,v)(w')=\frac{\ph(x+tw')-\ph(x)-t\la v,w'\ra}{\sm t^2}\ge-\rho\|w'\|^2.
\end{equation*}
This implies by passing to the limit as $w'\to w$ and $t\dn 0$ that \eqref{quad-est} holds. Since the function $w\mapsto\d^2\ph(x,v)(w)$ is positive homogeneous of degree $2$, we obtain $\d^2\ph(x,v)(0)=0$, which verifies that $\d^2\ph(x,v)$ is proper for every $(x,v)\in(\gph\sub\ph)\cap\B_{\ve/2}(\ox,\ov)$, and thus completes the proof.
\end{proof}\vspace*{-0.07in}

Next we are going to show that each subproblem \eqref{subp3} admits a unique solution under the tilt stability of a given local minimizer of the constrained optimization problem \eqref{coop}. \vspace*{-0.05in}

\begin{Proposition}[\bf solvability of subproblems in constrained optimization]\label{solv6} Let $\ph\colon\R^n\to\oR$ be taken from \eqref{coop2} with $\ox\in\dom\ph$, where $\ox$ is a tilt-stable local minimizer of $\ph$ with modulus $\kappa>0$. Suppose in addition that the assumptions of Proposition~{\rm\ref{tedi}} are satisfied with $\ov:=-\nabla\psi(\ox)$.
Then there exists $\ve>0$ such that for any $(x,v)\in(\gph\sub\ph)\cap\B_\ve(\ox,0)$ the unconstrained optimization problem
\begin{equation}\label{subp6}
\min\;\;\la v,w\ra+\sm\d^2\ph(x,v)(w)\;\;\mbox{subject to}\quad w\in \R^n
\end{equation}
admits a unique optimal solution.
\end{Proposition}\vspace*{-0.17in}
\begin{proof} It follows from \cite[Proposition~9.1]{mms19} that the metric subregularity assumption (a) in Proposition~\ref{tedi} yields the continuous prox-regularity of $\ph$ at $\ox$ for $0$. Employing Propositions~\ref{tedi} and \ref{prope} ensures the existence of $\ve>0$ such that for every $(x,v)\in(\gph\sub\ph)\cap\B_\ve(\ox,0)$ the function $\ph$ is twice epi-differentiable at $x$ for $v$ and that the second subderivative $\d^2\ph(x,v)$ is proper and lower semicontinuous. Proceeding as in the proof of Theorem~\ref{solv2} with the usage of \cite[Corollary~6.3]{pr96} and decreasing $\ve$ if necessary, we get that for any $(x,v)\in(\gph\sub\ph)\cap\B_\ve(\ox,0)$ the objective function in \eqref{subp6} is strongly convex with modulus $\sm\kappa^{-1}$. This surely verifies the existence of a unique optimal solution to subproblem \eqref{subp6} and hence completes the proof.
\end{proof}\vspace*{-0.07in}

To proceed with justifying superlinear convergence of Algorithm~\ref{nm3}, we need to investigate the semismooth$^*$ property of $\sub\ph$ for the objective function $\ph$ of \eqref{coop2}. This requires some additional assumptions. Recall  from \cite{bs} that a closed convex set $\Th\subset\R^m$ is ${\cal C}^2$-{\em cone reducible} at $\oz\in\Th$ to
a closed convex cone $Q\subset\R^s$ if there exist a neighborhood ${\cal O}\subset\R^m$ of $\oz$ and a ${\cal C}^2$-smooth mapping $h\colon\R^m\to\R^s$ such that
\begin{equation*}\label{red}
\Th\cap{\cal O}=\big\{z\in{\cal O}\big|\;h(z)\in Q\big\},\quad h(\oz)=0,\;\mbox{ and }\;\nabla h(\oz)\;\mbox{ has full rank}\;s.
\end{equation*}
It is proved in \cite[Theorem~6.2]{mms19} that ${\cal C}^2$-cone reducible sets are parabolically regular. Moreover, the latter result tells us that such sets satisfy assumption (b) in Proposition~\ref{tedi}.\vspace*{0.03in}

The next proposition reveals conditions on the initial data of \eqref{coop} ensuring the semismooth$^*$ property of the subgradient mapping $\partial\ph$ in \eqref{coop2}.\vspace*{-0.05in}

\begin{Proposition}[\bf propagation of semismooth$^*$ property in constrained optimization]\label{ssch} Taking the cost function $\ph$ in \eqref{coop2}, let $\ox\in\dom\ph$ with $0\in\sub\ph(\ox)$.
Assume that the convex set $\Th$ from \eqref{coop} is ${\cal C}^2$-cone reducible at $f(\ox)$, and that the nondegeneracy condition
\begin{equation}\label{nondeg}
\span\big\{N_\Th\big(f(\ox)\big)\big\}\cap\ker\nabla f(\ox)^*=\{0\}
\end{equation}
holds. If the normal cone mapping $N_\Th$ is semismooth$^*$ at $(f(\ox),\bar\lm)$, where $\bar\lm$ is the unique Lagrange multiplier in $\Lm(\ox,-\nabla\psi(\ox))$, then the subgradient mapping $\sub\ph$ is semismooth$^*$ at $(\ox,0)$.
\end{Proposition}\vspace*{-0.17in}
\begin{proof}
It follows from the elementary first-order subdifferential sum rule (see, e.g., \cite[Proposition~1.30(ii)]{m18}) that the stationary condition $0\in\sub\ph(\ox)$ amounts to $-\nabla\psi(\ox)\in N_\Th(\ox)$. It is also well known that the nondegeneracy condition \eqref{nondeg} implies that the Lagrange multiplier set $\Lm(\ox,-\nabla\psi(\ox))$ is a singleton; see, e.g., \cite[Proposition~4.75]{bs}. Employing the directional coderivative calculation from \cite[Theorem~4]{go17} tells us that for any $(u,v)\in T_{\scriptsize{\gph\sub\ph}}(\ox,0)$ we get
\begin{equation*}
\big(D^*\sub\ph\big)\big((\ox,0);(u,v)\big)(w)=\nabla_{xx}^2L(\ox,\bar\lm)w+\nabla f(\ox)^*D^*N_\Th\big((f(\ox),\bar\lm);(\nabla f(\ox)u,\xi)\big)\big(\nabla f(\ox)w\big)
\end{equation*}
whenever $w\in\R^n$, where $\xi\in\R^m$ is the unique solution to the system
\begin{equation*}
\xi\in DN_\Th\big(f(\ox),\bar\lm\big)\big(\nabla f(\ox)u\big),\;\;v=\nabla_{xx}^2L(\ox,\bar\lm)u+\nabla f(\ox)^*\xi.
\end{equation*}
Picking $(w,q)\in\gph(D^*\sub\ph)((\ox,0);(u,v))$ gives us $p\in D^*N_\Th((f(\ox),\bar\lm);(\nabla f(\ox)u,\xi))(\nabla f(\ox)w)$ such that $q=\nabla_{xx}^2L(\ox,\bar\lm)w+\nabla f(\ox)^*p$. Since $N_\Th$ is assumed to be semismooth$^*$ at $(f(\ox),\bar\lm)$, we deduce from \eqref{semi6} the representation
\begin{equation*}
\la\nabla f(\ox)w,\xi\ra=\la p,\nabla f(\ox)u\ra,
\end{equation*}
which in turn brings us to the equalities
\begin{equation*}
\la q,u\ra =\la\nabla_{xx}^2L(\ox,\bar\lm)w,u\ra+\la p,\nabla f(\ox)u\ra=\la\nabla_{xx}^2L(\ox,\bar\lm)w,u\ra+\la\nabla f(\ox)w,\xi\ra=\la w,v\ra.
\end{equation*}
According to  \eqref{semi6}, we arrive at the claimed semismooth$^*$ property of $\sub\ph$ at $(\ox,0)$.
\end{proof}\vspace*{-0.07in}

It is important to mention that the semismoothness$^*$ of $N_\Th$ imposed in Proposition~\ref{ssch} is automatically satisfied if $\Th$ is a polyhedral convex set, the second-order (Lorentz, ice-cream) cone, and the cone of positive semidefinite symmetric matrices. This results from the well-known fact that the projection mapping to such sets satisfies the estimate in Proposition~\ref{semi5}(ii) (see, e.g., \cite{css,s02}), which is equivalent to saying that it is semismooth$^*$. Using this  and the relationship $\Pi_\Th=(I+N_\Th)^{-1}$ for the projection operator, it is not hard to check that  the mapping $N_\Th$ is semismooth$^*$ at every point of its graph for the aforementioned convex sets.\vspace*{0.03in}

We are now in a position to present the main result of this section.\vspace*{-0.05in}

\begin{Theorem}[\bf superlinear convergence of the second subderivative-based Newton algorithm for constrained problems]\label{subd-semi} Let $\ox$ be a tilt-stable local minimizer of the cost function $\ph$ in \eqref{coop2}, and let the set $\Th$ in \eqref{coop} be ${\cal C}^2$-cone reducible at $f(\ox)$ under the fulfillment of the nondegeneracy condition \eqref{nondeg}. Assume also that the normal cone mapping $N_\Th$ is semismooth$^*$ at $(f(\ox),\bar\lm)$, where $\bar\lm$ is the unique Lagrange multiplier in $\Lm(\ox,-\nabla\psi(\ox))$. Then for any small $r>0$ there exists $\dd>0$ such that for each starting point $x_0\in\B_\dd(\ox)$ the {\rm(}unique{\rm)} sequence $\{x_k\}$ constructed by Algorithm~{\rm\ref{nm3}}
converges to $\ox$ and the rate of convergence is superlinear.
\end{Theorem}\vspace*{-0.17in}
\begin{proof} Since $\Th$ is ${\cal C}^2$-cone reducible at $f(\ox)$, the assumptions of Proposition~\ref{tedi}(b) hold by \cite[Theorem~6.2]{mms19}. Furthermore, the nondegeneracy condition \eqref{nondeg} yields the fulfillment of assumption (a) in Proposition~\ref{tedi}. Proposition~\ref{ssch} tells us that the subdifferential mapping $\sub\ph$ is semismooth$^*$ at $(\ox,0)$. It follows from Proposition~\ref{solv6} that  subproblem \eqref{subp3} admits a unique solution whenever $x$ is sufficiently close to $\ox$. As discussed above, the sequence $\{x_k\}$ generated by Algorithm~\ref{nm3} is actually induced by Algorithm~\ref{nm2} for the regularized function $e_r\ph$. This together with Theorem~\ref{supunc} implies that the uniquely determined sequence
$\{x_k\}$ constructed by Algorithm~\ref{nm3} converges to $\ox$ and the rate of convergence is superlinear.
\end{proof}\vspace*{-0.07in}

Observe that the choice of $v_k=\nabla(e_r\ph)(x_k)$ in the obtained conditions for solvability of subproblems and superlinear convergence of Algorithm~\ref{nm3} is the only one that is expressed not in terms of the given data of the constrained problem \eqref{coop} but via its Moreau regularization \eqref{moreau}. It is different from Algorithm~\ref{nm2} for the case of ${\cal C}^{1,1}$ functions, where $v_k=\nabla\ph(x_k)$. Although the relationships in \eqref{eq03} and \eqref{rd01} help us to write Algorithm~\ref{nm3} entirely in terms of the initial data of the constrained problem \eqref{coop}, the choice of $(x_k-rv_k,x_k)$ as a point from $\gph\partial\ph$ is a hard task numerically. We again refer the reader to the discussion at the end of Section~\ref{sect05} on the calculation of Moreau envelopes.\vspace*{0.03in}

Let us conclude this paper by some comments on differences between Algorithm~\ref{nm3} and the basic {\em sequential quadratic programming method} (SQP) to solve constrained optimization problems \eqref{coop}. The nondegeneracy condition \eqref{nondeg} implies that the Lagrange multiplier set $\Lambda(x_k-rv_k,v_k-\nabla\psi(x_k-rv_k))$, which appears in the calculation of $\d^2\ph(x_k-rv_k,v_k)(w)$ via \eqref{ssf}, is a singleton $\{\lm_k\}$. This, together with \eqref{ssf}, tells us that subproblem \eqref{subp3} in Algorithm~\ref{nm3} reduces to
\begin{equation}\label{subp6}
\underset{w\in\R^n}{\min}\;\la v_k,w\ra+\sm\big\la\nabla_{xx}^2L(x_k-rv_k,\lm_k)w,w\big\ra+\d^2\dd_\Th\big(f( x_k-rv_k),\lm_k\big)\big(\nabla f(x_k-rv_k)w\big).
\end{equation}

To compare Algorithm~\ref{nm3} with the basic SQP method, assume further that $\Th=\{0\}^{s}\times \R_-^{m-s}$, where $s$ is a positive integer with $0\le s\le m$. This choice of $\Th$ reduces the constrained problem \eqref{coop} to a nonlinear programming problem with $s$ equality constraints and $m-s$ inequality constraints. In this setting, it is well known that the nondegeneracy condition \eqref{nondeg} amounts to the classical linear independent constraint qualification (LICQ). Moreover, we know from \cite[Theorem~5.2]{mr} that the fulfillment of the LICQ implies that the tilt stability of a local minimizer is equivalent to Robinson's strong second-order sufficient condition. In summary, the nondegeneracy and tilt stability assumptions in Theorem~\ref{subd-semi} amount to the LICQ and the strong second-order sufficient condition for nonlinear programming problems, which were used conventionally for local convergence analysis of the basic SQP method for this class of problems. Note that the semismoothness$^*$ of $N_\Th$ for this choice of $\Th$ results from \cite[Proposition~3.5]{go19}, since $\gph N_\Th$ is a finite union of polyhedral convex sets. It is worth mentioning that the reason for us to impose the nondegeneracy condition in Theorem~\ref{subd-semi} is to ensure the semismoothness$^*$ of the constraint set $\O$ defined by \eqref{coop2}. Having a chain rule for the semismoothness$^*$ of the latter set under weaker constraint qualifications would allow to improve Theorem~\ref{subd-semi}.

Now let us compare subproblem \eqref{subp6} with that of the basic SQP method. To this end, we begin with the simplification of \eqref{subp6}. It follows from \cite[Example~3.4]{mms19} that
\begin{equation}\label{ssp}
\d^2\dd_\Th\big(f(x_k-rv_k),\lm_k\big)\big(\nabla f(x_k-rv_k)w\big)=\dd_{K_\Th(f(x_k-rv_k),\lm_k)}\big(\nabla f(x_k-rv_k)w\big).
\end{equation}
To obtain a convenient formula for the critical cone of $\Th$, pick $(z,\lm)\in\gph N_\Th$. Since $\Th=\{0\}^{s}\times \R_-^{m-s}$, it gives us $\lm\in\R^{s}\times\R^{m-s}_+$.
Define further the index sets
\begin{equation*}
I (z)=\big\{i\in\{s+1,\ldots,m\}|\;z_i=0\big\}\quad \mbox{and}\quad I_+(z,\lm )=\big\{i\in I(z)|\,\lm_i>0\big\},
\end{equation*}
where $z=(z_1,\ldots,z_m)$ and $\lm=(\lm_1,\ldots,\lm_m)$. Using these index sets, we can conclude that
\begin{equation*}
(y_1,\ldots,y_m)\in K_\Th (z,\lm)=T_\Th(z)\cap\{\lm\}^\perp\iff
\begin{cases}
y_i=0&\mbox{for all}\;\;i\in\{1,\ldots,s\}\cup I_+(z,\lm),\\
y_i\le 0&\mbox{for all}\;\;i\in I(z)\setminus I_+(z,\lm).
\end{cases}
\end{equation*}
Combining this representation with \eqref{ssp} allows us to equivalently express subproblem \eqref{subp6} as
\begin{equation*}
\begin{cases}
\underset{w\in\R^n}{\min}&\la v_k,w\ra+\sm\big\la\nabla_{xx}^2L(x_k-rv_k,\lm_k)w,w\big\ra\\
\mbox{subject to}&\nabla f_i(x_k-rv_k)w=0,\;\;i\in\{1,\ldots,s\}\cup I_+(z_k,\lm_k),\\
&\nabla f_i(x_k-rv_k)w\le 0,\;\;i\in I(z_k)\setminus I_+(z_k,\lm_k ),
\end{cases}
\end{equation*}
where $f=(f_1,\ldots,f_m)$ with $f_i:\R^n\to\R$ and $z_k:=f(x_k-rv_k)$. In contrast, subproblems of the basic SQP method in this setting at the current primal-dual iterate $(x_k,\lm_k)$ can be formulated as follows (see, e.g., \cite[Section~6.3]{dr}):
\begin{equation*}
\begin{cases}
\underset{w\in\R^n}{\min}&\la\nabla\psi(x_k),w\ra+\sm\big\la\nabla_{xx}^2L(x_k,\lm_k)w,w\big\ra\\
\mbox{subject to }&f_i(x_k)+\nabla f_i(x_k)w=0,\;\;i=1,\ldots,s,\\
&f_i(x_k)+\nabla f_i(x_k)w\le 0,\;\;i=s+1,\ldots,m.
\end{cases}
\end{equation*}
We see that the main difference between the subproblems in Algorithm~\ref{nm3} and in the basic SQP method is that the latter keeps the constraint set $\Th$ {\em untouched}, while our algorithm replaces it by the critical cone associated with each iteration of Algorithm~\ref{nm3}. Finally, note that while SQP methods generate a primal-dual sequence, Algorithm~\ref{nm3} constructs only a primal sequence. In this regard, these algorithms behave differently.\\[1ex]
{\bf Acknowledgements}. The authors are very grateful to two anonymous referees and the handling Associate Editor for their constructive remarks and suggestions, which allowed us to significantly improve the original presentation.
\end{document}